\documentclass[ijaa]{informs4}
\usepackage{eqndefns-left} 

\usepackage[T1]{fontenc}      
\usepackage[utf8]{inputenc}   
\usepackage{textcomp}         

\usepackage{hyperref}
\hypersetup{hidelinks}

\RequirePackage{tgtermes}
\RequirePackage{newtxtext}
\RequirePackage{newtxmath}
\RequirePackage{bm}
\RequirePackage{endnotes}
\usepackage{listings}

\usepackage{float}
\newfloat{codefloat}{h}{loc}
\floatname{codefloat}{Code}

\usepackage{caption}
\usepackage[font=small, labelfont=normalfont, textfont=normalfont, skip=10pt]{caption}

\makeatletter
\renewcommand{\figurename}{Figure}               
\renewcommand{\fnum@figure}{\figurename~\thefigure} 
\makeatother

\OneAndAHalfSpacedXII 

\usepackage{algorithm}
\usepackage{algpseudocode}
\usepackage{tikz}
\usepackage{booktabs} 

\usepackage{natbib}
 \bibpunct[, ]{(}{)}{,}{a}{}{,}%
 \def\bibfont{\small}%
 \def\BIBand{and}%

\EquationsNumberedThrough    

\TheoremsNumberedThrough     
\ECRepeatTheorems  %

\MANUSCRIPTNO{IJAA-0001-2024.00}

\def\theARTICLETOPLEFT{}
\def\theARTICLETOPRIGHT{}

\begin{document}


\RUNAUTHOR{Llenas et al.}

\RUNTITLE{PepsiCo Deploys AI-Driven Pricing and Promotion Optimization at Scale}

\TITLE{PepsiCo Deploys AI-Driven Pricing and Promotion Optimization at Scale}


\ARTICLEAUTHORS{%
\AUTHOR{Aleix Llenas \textsuperscript{a,*}, Eduardo Salazar-Treviño \textsuperscript{a}, Francisco Leskovar \textsuperscript{a}, Francesc Pons-Llopis \textsuperscript{a}, Federico Todeschini \textsuperscript{a}, Deepika Gowda \textsuperscript{a}, David Bofill \textsuperscript{a}, Laxmi Anish  \textsuperscript{a} and Michael Cleavinger \textsuperscript{a,**}}
\vspace{5pt}
\AFF{\fs.10.12. \textsuperscript{a} Data Science \& Advanced Analytics, PepsiCo, Inc.\\
\vspace{5pt}
\textsuperscript{*}\EMAIL{aleix.llenasfarras@pepsico.com}; 
\textsuperscript{**}\EMAIL{michael.cleavinger@pepsico.com}}
\vspace{15pt}
\AFF{{\fs.10.12. Published in {\bi\theJOURNAL}.}\\
\vspace{6pt}
{\fs.10.12.Accepted May 18, 2026 $\cdot$ DOI: \href{https://doi.org/10.1287/inte.2025.0302}{\textcolor{blue}{10.1287/inte.2025.0302}}}}
}

\ABSTRACT{%
Effective pricing and promotion planning constitutes a central pillar of strategic revenue management for firms operating in highly competitive and dynamic markets. These planning activities require the simultaneous consideration of demand elasticity, competitor actions, channel and market specific constraints, and financial objectives. As the dimensionality and interdependencies inherent in these problems increase, manual or traditional approaches become suboptimal and insufficient. In this context, Operations Research provides a robust methodological foundation for scalable, data-driven decision support systems that can optimize complex planning processes across large product and customer portfolios.

This paper presents two large-scale optimization systems developed and deployed at PepsiCo to support Revenue Growth Management initiatives: PromoAI and PricingAI. PromoAI integrates machine learning-based promotional forecasts with a mixed-integer linear programming model to optimize promotional calendars across trade channels. The system navigates through millions of product-promotion-timing combinations to find the one that maximizes PepsiCo and retailer revenues, subject to a wide range of customizable business constraints encoded in a modular, user-configurable interface.

On the other hand, PricingAI focuses on the optimization of base prices across product portfolios over multi-period horizons. The system employs Bayesian hierarchical models to estimate own and cross-price elasticities and captures competitive interactions at the product level. These elasticity estimates are fed into a nonlinear programming optimization engine that recommends price changes aligned with revenue and margin targets, while incorporating operational constraints such as price thresholds, volume or profit margins, and channel- and market-specific business rules.

Together, these systems demonstrate the feasibility and scalability of advanced optimization in large-scale enterprise environments. They highlight the value of integrating statistical learning with mathematical programming to enable enterprise-level, automated decision-making that is both data-informed and aligned with strategic business objectives.

}%




\maketitle



\section{Introduction}\label{sec:Intro}

Consumer packaged goods (CPG) companies operate in highly competitive and margin-sensitive environments, where strategic decisions related to pricing and promotional planning play a critical role in driving revenue, profitability, and market share. PepsiCo, Inc., a global leader in the CPG sector with annual revenues exceeding \$90 billion, manages a vast portfolio of products across numerous geographies and retail channels \citep{PepsiCo}. Within such a complex landscape, decisions surrounding price setting and trade promotions must account for heterogeneous consumer demand, retailer dynamics, seasonal trends, and competitor actions. 

In the CPG industry, revenue management relies heavily on two distinct but interrelated levers: trade promotion planning and strategic base pricing. While both aim to maximize revenue and profitability, they operate on different timelines and involve fundamentally different decision structures.

Trade promotions are temporary price reductions or special offers (e.g., "Buy One Get One," "2 for \$5") agreed upon between the manufacturer and the retailer. The core business problem is the construction of an optimal promotional calendar for a specific retailer over a defined horizon, typically a quarter or a year. This is not merely a question of which products to discount, but a complex scheduling problem subject to strict operational constraints. Retailers have limited "slotting" capacity for promotions—for example, a supermarket chain may only allow a specific brand category to appear on the front page of a circular or on an end-of-aisle display for a limited number of weeks per year. This calendar serves as a formal agreement between PepsiCo and the retailer to align their commercial efforts. In the case a retailer is unwilling or unable to execute a specific promotion, they may opt-out, provided this is managed through a formal negotiation.

Planners must determine the optimal timing, depth, and mechanics of these events for hundreds of Product Promotion Groups (PPGs). A valid calendar must respect exclusivity rules (e.g., avoiding simultaneous promotions of conflicting pack sizes), adhere to minimum spacing requirements between events (to prevent "pantry loading" where consumers stock up and delay future purchases), and align with major seasonal events like the Super Bowl or holidays. The objective is to maximize the "lift"—the incremental volume generated above baseline sales—while ensuring the cost of the promotion (trade spend) does not erode profit margins below acceptable thresholds for either PepsiCo or the retail partner.

In contrast to the tactical nature of promotions, base pricing involves setting the standard "everyday" price of products. These decisions are strategic, typically occurring on a semi-annual or annual cycle, and carry significant long-term risk. The fundamental challenge here is estimating and managing price elasticity—the sensitivity of consumer demand to permanent price changes.

The decision-maker must define a price list for the upcoming year that balances volume retention with margin expansion. This problem is complicated by the portfolio effect: increasing the price of a premium brand may drive consumers to a value brand within the same portfolio (cannibalization) or to a competitor’s product (switching). Furthermore, pricing in CPG is rarely continuous; it is often constrained by "price ladders" (logical gaps between different pack sizes) and psychological price thresholds (e.g., keeping a single-serve item under \$1.99 or 20 Pesos). The business problem is to find a set of prices across the entire portfolio that achieves corporate financial targets without triggering a price war or ceding market share to competitors.

Traditionally, CPG companies have relied on manual processes, domain expertise, and historical heuristics to inform pricing and promotion planning \citep{ma2021dynamic, digital_cpg}. These methods reflect deep business intuition, accumulated market knowledge, and a strong understanding of local dynamics. However, while such approaches are valuable and often effective in specific contexts, they may struggle to scale efficiently or capture complex interactions across a growing number of variables. Another common challenge is the lack of a unified methodological framework across global markets: different countries or business units often employ entirely distinct approaches to revenue planning, limiting consistency, comparability, and the potential for centralized optimization.

Recent advances in Artificial Intelligence (AI), Machine Learning (ML), and Operations Research (OR) present significant opportunities to enhance revenue management practices \citep{bertsimas2020predictive}. ML models can improve demand forecasting, quantify price and promotion elasticities, and capture competitive effects with greater accuracy. OR techniques, particularly in mathematical optimization, offer structured methods to generate prescriptive recommendations that balance business objectives and business rules in a transparent and repeatable way. Together, these tools provide a powerful extension to traditional methods, enabling companies to move toward more scalable, data-informed, and globally harmonized decision-making processes.

This study presents two AI-driven initiatives for Revenue Growth Management (RGM) developed within PepsiCo's Data Science department: PromoAI and PricingAI. PromoAI is an advanced optimization platform designed to automate and enhance trade promotion planning through the generation of promotional calendars. Its goal is to identify promotion schedules that maximize financial performance for both PepsiCo and its retail partners. Depending on the market, business units design promotional calendars on a quarterly, semiannual, or annual basis for execution at the retailer level.

In parallel, PricingAI enables granular, data-driven pricing decisions by estimating price elasticities across multiple dimensions, such as brand, package size, and retailer. Its primary objective is to develop optimized base pricing strategies that drive both revenue growth and profitability. Adjustments to base product prices also vary by market, but typically occur on a semiannual or annual cycle.

Revenue management has been extensively studied across the fields of operations research and machine learning, with a broad range of techniques proposed for price optimization, promotional planning, and product assortment design \citep{talluri2006theory, fattahi2022customer}. The INFORMS Journal on Applied Analytics (formerly Interfaces) has a rich history of documenting such applications, including the group pricing optimization system at Marriott \citep{hormby2010marriott} and the integrated demand and price optimization framework at Alibaba \citep{deng2023alibaba}. 

Other modern approaches include reinforcement learning for dynamic pricing in online retail environments \citep{ferreira2018online, cohen2018dynamic}, prescriptive optimization frameworks that embed machine learning-based demand forecasts into mathematical programming formulations \citep{ito2017optimization}, and choice-based models for joint assortment and pricing under substitution behavior \citep{miao2020data}. While these methods have shown strong theoretical and empirical performance in controlled or small-scale settings, relatively few studies have demonstrated their effectiveness in large-scale, enterprise environments.

This work addresses that gap by presenting two enterprise-grade optimization systems developed and deployed at PepsiCo. Both systems are built to operate at scale—across product portfolios, geographies, and planning cycles—and are designed with explicit attention to governance, human-in-the-loop workflows, and organizational adoption as first-class design goals alongside technical performance.

A key contribution of this work lies in demonstrating how the thoughtful integration of AI and OR methods can effectively augment traditional revenue management practices within a complex, global CPG organization. We describe the architectural and modeling approximations necessary to ensure computational tractability, the alignment of model outputs with business decision frameworks, and the scalability of these tools under heterogeneous market constraints. The aim of this paper is to provide a documented case study of how advanced optimization and machine learning techniques can be applied and sustained in practice — delivering measurable business value within a Fortune 100 company.

\section{Methodology}\label{sec:Method}

The optimization systems developed at PepsiCo employ methodological approaches tailored to the distinct requirements of promotional planning and pricing decisions. Despite these differences, both systems share a common architectural framework: ML-based demand estimation followed by a mathematical optimization stage. 

The systems diverge, however, in their technical implementations and optimization strategies. PromoAI formulates promotional planning as a constrained calendar optimization problem solved via mixed-integer linear programming (MILP), whereas PricingAI adopts a metaheuristic approach, using differential evolution (DE) to identify optimal pricing strategies.

MILP is a mathematical optimization technique that extends linear programming by allowing decision variables to be constrained to integer values, either partially (mixed-integer) or fully \citep{wolsey1999integer}. MILP models consist of a linear objective function subject to linear equality and inequality constraints, with some variables required to take binary or integer values. This framework is particularly well-suited for planning problems that involve discrete choices, such as whether to run a promotion in a specific week or whether to include a product in a promotional calendar. The strength of MILP lies in its ability to capture complex business logic (e.g., exclusivity rules, minimum frequency of promotions, budget limits) within a rigorous and solvable optimization structure. Modern MILP solvers (e.g., Gurobi, CPLEX) have become increasingly efficient, capable of handling large-scale, real-world problems with thousands of variables and constraints \citep{bertsimas1997introduction}. In the context of revenue management, MILP provides a robust approach for prescriptive planning when the underlying relationships can be well approximated with linear formulations \citep{wicaksono2008piecewise}.

DE is a population-based metaheuristic optimization algorithm particularly well-suited for continuous global optimization problems with non-convex, non-differentiable objective functions \citep{storn1997differential}. Unlike gradient-based methods that can become trapped in local optima, DE maintains a population of candidate solutions that evolve through mutation, crossover, and selection operations. The algorithm's key strength lies in its simple yet powerful mutation strategy, which creates new candidate solutions by adding weighted differences between population members to a base vector. This self-organizing property enables DE to adapt its search behavior to the landscape of the objective function, making it remarkably effective for complex optimization problems. In revenue management contexts, DE is particularly valuable for pricing optimization where the objective function—incorporating elasticity models, competitive dynamics, and business constraints—exhibits high nonlinearity and multiple local optima \citep{talluri2006theory}. The algorithm's parameter efficiency, the natural handling of constraints and the robust performance across diverse problem types have made it a good candidate for real-world applications where traditional optimization methods struggle.

In addition to the optimization methods discussed, each one works in combination with an ML model operated on the previous stages. 

Before the MILP solver we have trained an ML forecast model, usually of the Gradient Boosting Machines (GBMs) class. GBMs are ensemble learning methods that construct predictive models by sequentially combining weak learners, typically decision trees, to minimize a specified loss function. Each tree is trained to correct the residual errors of the ensemble formed by its predecessors, allowing the model to capture complex, nonlinear relationships and high-order feature interactions. GBMs are widely used for structured data forecasting tasks due to their strong predictive performance, robustness to outliers, and ability to handle mixed data types without extensive preprocessing \citep{friedman2001greedy}. Popular implementations include XGBoost, LightGBM, and CatBoost, which offer efficient training and built-in regularization techniques to mitigate overfitting.

On the other hand, on top of the DE algorithm we have executed a Bayesian ML model. This is a probabilistic modeling framework that treats model parameters as random variables and uses Bayes' theorem to update their distributions in light of observed data. Unlike point-estimate approaches, Bayesian methods produce full posterior distributions, enabling explicit representation of uncertainty and incorporation of prior knowledge \citep{gelman1995bayesian}. This is particularly valuable in applications with limited data or high variance, where capturing the confidence of predictions is as important as the predictions themselves. Bayesian regression and hierarchical models are commonly used to infer relationships between variables while accounting for latent structure, making them well-suited to domains where interpretability and uncertainty quantification are critical.

The choice of distinct optimization paradigms—MILP for promotions versus Bayesian inference with metaheuristic optimization for pricing—reflects fundamental differences in these decision domains. Promotional planning operates in a data-rich environment with discrete, reversible decisions and well-defined constraints. In contrast, strategic pricing faces severe data limitations due to infrequent price changes, requires explicit uncertainty quantification given the risks of incorrect pricing, and must navigate highly non-convex objective functions arising from complex elasticity interactions. The following sections present an in-depth analysis of the technical characteristics of each project.

\subsection{PromoAI: Promotional Planning Optimization}
\subsubsection{Data Preparation}
Before the development of any machine learning or optimization models, a fundamental prerequisite for this project is the availability of comprehensive historical datasets from the retailer. Specifically, models require detailed records of weekly sales, pricing, and promotional activities. In addition, information on product placement—such as shelf location (upper, middle, or lower)—is essential for some markets. Financial data, including trade spending and product cost, is also necessary to compute margins and profitability. 

Collecting and preparing these data presented several challenges. Retailers differ in how they define promotional events, record price changes, and structure product hierarchies, which requires substantial harmonization across markets. For example, some retailers run promotions at the individual product level, whereas others apply a single promotional action simultaneously to groups of products.

Historical records often contained missing weeks, inconsistent promotion labels, and price entries that did not reflect effective in-store prices due to delayed system updates. In addition, product identifiers changed over time due to new product introductions or Stock Keeping Unit (SKU) rationalization, requiring the creation of consistent longitudinal product mappings.

The data preparation process therefore included extensive cleaning and validation steps. These included imputing missing values, correcting pricing anomalies and validating sales spikes against documented promotional activity. Automated data quality checks were implemented to flag implausible price–volume combinations and margin inconsistencies. Only after these harmonization and validation steps were completed did the machine learning modeling phase begin.
\subsubsection{Forecasting Model}

PromoAI employs a global regression model using LightGBM \citep{ke2017lightgbm} to predict promotional demand. Rather than building separate models for each product, we adopt a unified approach where all Product Promotion Groups (PPGs) are combined in a single dataset. A PPG represents a group of SKUs with similar characteristics, allowing the model to share information across related products while still learning distinct baseline demand levels for each group.

The model learns demand as a function of several key drivers, including:
\begin{itemize}
    \item \textbf{Temporal variables}, such as week, month, and holiday indicators,
    \item \textbf{Promotional characteristics}, including discount depth, type of promotional
          mechanism (e.g., percent discount, multi-buy offers), and product placement,
    \item \textbf{Product attributes}, such as brand and category.
\end{itemize}

By modeling all PPGs jointly, PromoAI captures common demand behaviors across the portfolio and enables more robust demand predictions, even for products with limited historical data.

The LightGBM model is trained using a mean absolute percentage error (MAPE) objective, which emphasizes relative forecast accuracy and is well-suited to the heterogeneous volume scales observed across different PPGs and markets. Hyperparameter selection is performed via randomized cross-validation over a held-out temporal split (the most recent weeks are reserved as validation data to mimic the forward-looking nature of planning). The table describing the hyperparameters selected for this specific instance of the model can be found in the Appendix A.1. 

The learned model structure can be partially characterized through feature importance scores derived from the ensemble. The most influential feature is the PPG identifier, which encodes baseline demand level and brand strength, reflecting that intrinsic product characteristics are the dominant driver of demand variation. Temporal features rank prominently, with year, holiday indicators, week-of-year, and quarter all appearing in the top ten, capturing both long-term trend and seasonal demand patterns. Promotional characteristics complete the ranking, including mechanic type, retailer base price, actual transaction price, promotional depth, and month, reflecting the direct price stimulus and mechanic-specific response patterns that drive promotional lift.

\subsubsection{Piecewise Linear Demand Approximation and Breakpoint Selection} \label{pwl_section}

Before describing the business constraints, it is important to elaborate on how the nonlinear demand forecasts produced by the LightGBM model are embedded into the MILP framework, since this step is both technically central and practically challenging.

The LightGBM model produces a nonlinear mapping from the feature vector to a demand forecast. Since MILP requires linear relationships, this nonlinear function is approximated via piecewise linear (PWL) segments in the discount pressure variable. Specifically, for each PPG, week, and promotion, the system evaluates the LightGBM model at a grid of discount pressure values to construct the breakpoint vectors in Equation \ref{eq:pwl} of Appendix A. Figure~\ref{fig:pwl} illustrates this approximation for a representative PPG demand curve, showing how
three optimally placed breakpoints closely follow the nonlinear LightGBM response across the full range of competitive discount pressure values.

\begin{figure}[h]
\centering
\includegraphics[width=0.5\textwidth]{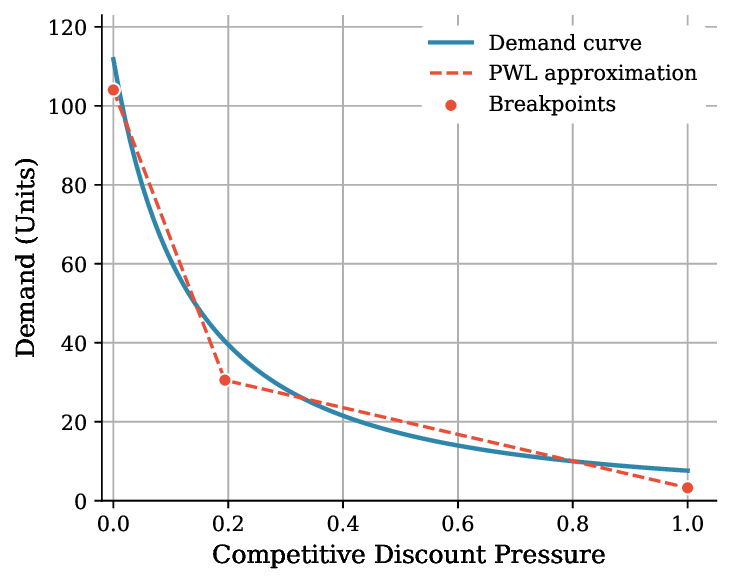}
\caption{Example of a piecewise linear (PWL) approximation of a representative demand response curve as a function of competitive discount pressure. Two segments with optimally placed breakpoints closely approximate the nonlinear demand response, enabling the MILP solver to consider the demand behavior while maintaining a tractable linear formulation.}
\label{fig:pwl}
\end{figure}

Breakpoint positions and the number of segments are determined through a two-stage automated optimization process. First, for a given number of breakpoints, the x and y coordinates of the breakpoints are jointly optimized using constrained nonlinear programming (SLSQP) to minimize the mean squared error between the piecewise linear approximation and the full LightGBM response curve. Second, this optimization is repeated across a range of segment counts, and the optimal number of breakpoints is selected using a knee-point detection algorithm that identifies the point of diminishing returns on the approximation error curve — the segment count beyond which additional breakpoints yield negligible improvement in fit relative to the additional solver complexity they introduce. In practice, 2--4 breakpoints have proven sufficient across all deployed markets, reflecting the relatively smooth and monotonic nature of the demand-discount pressure relationship in CPG retail, as illustrated in Figure~\ref{fig:pwl}. Once validated for a given market, the breakpoint configuration is retained across planning cycles unless a significant change in demand behavior warrants re-evaluation.

\subsubsection{Business Constraints}

The optimization engine schedules promotions over a multi-week planning horizon. For each product and each week, the model decides which promotional option (including ``no promotion'') should be selected.

As discussed before, the optimization problem is formulated as an MILP problem. The optimization objective balances PepsiCo and retailer goals. Users choose the weights that prioritize revenue versus margin, and PepsiCo versus retailer outcomes. For example, they may choose to maximize only PepsiCo revenue or select a balanced configuration that jointly improves PepsiCo and retailer performance. These weights are configured directly in the user interface, and the system automatically translates them into the optimization objective. Further details can be found in Appendix A.

PromoAI incorporates two distinct types of constraints: structural modeling constraints that are required to correctly formulate the calendar optimization problem, and business constraints that encode market-specific commercial rules configured by the user.

The structural constraints are common to all deployments and define the mathematical backbone of the problem:

\begin{enumerate}
    \item \textbf{Promotion exclusivity per week.}
          Each product must be assigned exactly one promotional option for each week in the planning
          horizon, ensuring a complete and operationally valid promotional calendar.

    \item \textbf{Competitive discount pressure.}
          To represent the influence of competing products, PromoAI computes an internal measure of
          ``discount pressure,'' which reflects the average promotional aggressiveness of products in
          similar size or format segments. The optimization model uses this measure to ensure that
          demand forecasts correctly reflect market-wide promotional intensity.

    \item \textbf{Demand approximation through piecewise linearization.}
          This is explained in the previous section \ref{pwl_section}.
\end{enumerate}

Beyond these structural constraints, PromoAI's flexibility is enabled by a modular library of business constraint templates, combined dynamically depending on the market configuration selected in the UI (User Interface). These fall into three categories:

\begin{itemize}
    \item \textbf{Financial constraints:} revenue targets, trade-spend limits, margin preservation,
          and market-share requirements.

    \item \textbf{Calendar constraints:} seasonal rules, holiday alignment, minimum and maximum
          promotional frequency requirements, and spacing rules that govern how often promotions can
          occur.

    \item \textbf{Execution constraints:} competitor lockout rules, ad-block linking, retail price
          restrictions, and minimum front-page exposure requirements.
\end{itemize}

Users configure these rules through the system interface. The rules are encoded in a JavaScript Object Notation (JSON) file, which PromoAI converts into the corresponding mathematical constraints

\subsection{PricingAI: Strategic Price Optimization}

Unlike promotional planning's discrete decision space, pricing optimization confronts fundamental uncertainties in demand response estimation. Price changes are infrequent in CPG markets, creating data scarcity that traditional econometric methods struggle to address. PricingAI therefore employs a two-stage approach: first, a Bayesian hierarchical model quantifies own- and cross-price elasticities; second, a metaheuristic optimization algorithm navigates the resulting non-convex landscape to identify robust pricing strategies.

\subsubsection{Data Preparation} 

PricingAI draws on three input data sources orchestrated through a product master file: historical weekly sales data (observed prices, sales, volumes, promotional activity, and distribution metrics at SKU-banner-week level), financial data (costs and sell-in prices for margin calculations), and conjoint survey data (when available, providing prior estimates for the elasticity components). Sell-in prices are the prices at which PepsiCo sells to retailers, as distinct from sell-out or shelf prices paid by consumers. The product master maps each SKU to its model characteristics, including brand, sub-brand, size, taste, and price line, and serves as the linking key across all input sources.

The main feature-engineering steps are the derivation of a regular price and the flagging of promotional weeks. Regular price is estimated as the maximum observed price within an 8-week sliding window (backward and forward looking), and a week is flagged as promotional if the observed price falls more than 5\% below this reference. Data is then aggregated from store-region level to a semi-regional level. When aggregating promotional flags, a revenue-weighted threshold of 30\% is applied: a SKU-banner-week is classified as promotional if more than 30\% of regional revenue was generated under promotional conditions in that week.

\subsubsection{Elasticity Estimation via Bayesian Hierarchical Model}

For each SKU/semiregion/week tuple, PricingAI assumes a multiplicative relationship between volume and its drivers. The model is estimated using a Bayesian hierarchical framework, which generates full posterior distributions over all parameters rather than single-point estimates. This enables principled uncertainty quantification over own- and cross-price elasticity estimates across the full product portfolio. The posterior mean of each elasticity component is then used as input to the downstream pricing optimization.

CPG retail data presents two structural estimation challenges: long periods of price stability create limited variation in the primary regressor, and the number of elasticity parameters to estimate across SKUs and banners vastly exceeds the number of independent price observations available for any individual product.

To address this, the model assumes a hierarchical relationship for price elasticity, decomposing own-price sensitivity across four shared attribute dimensions: brand, package size, flavor, and banner (retailer). For own-price elasticity these components are additive, while for cross-price elasticity they are multiplicative and bounded between 0 and 1. This structure reduces the effective parameter count from one free coefficient per SKU-banner pair to a compact set of shared components estimated jointly across the full portfolio, improving identification and generalization. It also reflects the economic intuition that similar products sold at the same retailer tend to exhibit similar demand responses. Table~\ref{tab:elast_components} summarizes the decomposition components.

\begin{table}[h]
\centering
\caption{Components of the own-price elasticity decomposition in PricingAI. 
All components are additive.}
\label{tab:elast_components}
\begin{tabular}{p{4cm}p{3.5cm}p{5.5cm}}
\hline
\textbf{Component} & \textbf{Scope} & \textbf{Role} \\
\hline
Grand elasticity      & 1 per category   & Category-level baseline \\
Brand elasticity      & 1 per brand      & Brand-family price sensitivity \\
Size elasticity       & 1 per pack size  & Pack-size sensitivity \\
Taste elasticity      & 1 per flavor     & Flavor-level sensitivity \\
Banner elasticity     & 1 per banner     & Retailer/channel sensitivity \\
\hline
\end{tabular}
\end{table}

Cross-price elasticities are also decomposed multiplicatively across sub-brand, size, and taste dimensions, reducing the parameter space from one free coefficient per SKU pair to a small set of attribute-level terms while preserving the relevant substitution patterns. 

Beyond own-price effects, the model captures the promotional discount elasticity, a binary promotional flag effect, distribution, monthly seasonality and cross-price elasticity terms. Additionally, the model includes both recency (more weight for more recent event) and volume importance (more weight to bigger SKUs).

When conjoint survey data is available, it provides informative prior means for each elasticity component. Tighter prior standard deviations on brand and taste dimensions reflect greater confidence in conjoint estimates along those attributes, while wider values on size and banner give more weight to the historical data. A softmax penalty enforces sign consistency by strongly discouraging positive own-price elasticity values across all SKUs regardless of data sparsity. 

The model is implemented in STAN \citep{carpenter2017stan}. Given the high dimensionality of the data, the model uses Automatic Differentiation Variational Inference (ADVI) \citep{kucukelbir2017advi}, maintaining a practical balance between statistical rigor and computational efficiency (in the range of hours). 

Model validation follows a three-stage protocol: (i) in-sample fit metrics at the SKU-banner-week level, with acceptance thresholds of $R^2 > 0.6$, wMAPE $< 0.4$, and $|\text{bias}| < 10\%$; (ii) business validation of elasticity magnitudes at brand-unit/size aggregation, with own-price elasticities expected in the range $(-2.5,\, 0)$ and positive cross-price elasticities consistent with substitution behavior; and (iii) end-to-end testing of the full optimization pipeline on the validated elasticity matrix. These validated elasticity matrices serve as the direct inputs to the price optimization engine described below.

\subsubsection{Price Optimization via Differential Evolution}

Given the non-convex nature of the pricing problem, we employ a DE method \citep{storn1997differential}, a population-based metaheuristic well-suited for global optimization. The procedure is implemented using algorithms from the SciPy optimization library \citep{2020SciPy-NMeth}, which offers a robust framework for large-scale numerical optimization. 
The core challenge comes from the objective function structure.

Similar to PromoAI, PricingAI builds its optimization objective based on the configuration selected by the user in the interface. At the core of this objective is a prediction of how demand changes when prices change, using elasticity estimates derived from the Bayesian model. These elasticity relationships are inherently nonlinear and involve interactions that create a highly irregular, non-convex optimization landscape.

PricingAI evaluates how different price configurations affect expected sales volume and then combines these volume predictions with margin calculations to estimate financial outcomes. Users specify whether the objective should prioritize revenue, profit, or a combination of the two, and the system automatically adjusts the objective accordingly before running the optimization.

PricingAI implements two primary constraint categories configured through the UI:

\begin{enumerate}
\item \textbf{Price Ladders}: Maintains logical price relationships within brand families of different sizes.
\item \textbf{Financial Constraints}: 
   \begin{itemize}
   \item Revenue targets
   \item Margin preservation
   \item Market share
   \item Profit constraint
   \item Bounds on volume (number of units sold)
   \end{itemize}
All these constraints can be set with a lower and upper bound that the user can configure through a dedicated User Interface (UI). See Appendix C for further details. 
\end{enumerate}

\section{System Architecture and Deployment}\label{sec:deployment}

Both PromoAI and PricingAI are implemented as modular Python packages and deployed on a shared Azure cloud infrastructure. The common operational workflow is initiated by business users interacting with a web-based UI, where they configure scenarios and submit optimization requests. These requests, converted to JSON format, are managed and orchestrated through Azure Data Factory pipelines, which connect the UI to the backend compute environment residing on Databricks notebooks. The notebooks execute the required optimization workflows by accessing underlying data for model training and building.

A crucial architectural principle for both systems is decoupling, which is fundamental to enhancing scalability, maintainability, and independent evolution.

\begin{itemize}
    \item \textbf{Loose Coupling}: The architecture maintains loose coupling between the UI, optimization engine, and data layers. This enables independent scaling and technological upgrades; for instance, the UI can be modernized, or the data pipeline upgraded, without affecting the core optimization algorithms.
    \item \textbf{Separation of Concerns}: ML models, used for inputs like demand predictions, are trained and versioned separately from the optimization engines. This allows data science teams to rapidly iterate on prediction models without impacting the stability of the optimization code.
    \item \textbf{Externalized Business Logic}: For both projects, business rules and market-specific constraints are externalized in JSON-based configuration files rather than being embedded in the code (an example can be seen in Appendix B). This configuration approach dramatically reduces the technical effort required for market expansion. A new market is onboarded by creating a configuration file that instantiates existing constraint modules with market-specific parameters. This design empowers business users to modify business rules without developer intervention, just by modifying the JSON file, ensuring the system can evolve rapidly with changing business needs.
\end{itemize}

\subsection{PromoAI}

PromoAI has been deployed across dozens of markets spanning multiple continents, including Canada, the United States, Mexico, Brazil, the United Kingdom, Australia, and South Africa. The code flexibility enables it to handle widely varying market characteristics and business requirements. Table~\ref{tab:promoai_scale} summarizes the range of problem sizes encountered across different markets.

\begin{table}[htbp]
\centering
\caption{PromoAI Problem Characteristics Across Markets}
\label{tab:promoai_scale}
\begin{tabular}{lcc}
\hline
\textbf{Characteristic} & \textbf{Minimum} & \textbf{Maximum} \\
\hline
Number of PPGs & 6-8 & 90+ \\
Planning Horizon & 8 weeks & 52 weeks \\
Decision Variables & \multicolumn{2}{c}{Hundreds of thousands} \\
Constraints & \multicolumn{2}{c}{Hundreds of thousands} \\
\hline
\end{tabular}
\end{table}

The significant variation in problem sizes reflects the diversity of market structures and business needs. Smaller markets with 6-8 PPGs typically represent focused product portfolios or specific retail channels, while larger markets encompass comprehensive product ranges across multiple categories.

The optimization problem is executed on a dedicated server running Gurobi version 12 \citep{gurobi}. Solution times vary substantially with problem size and complexity, ranging from a few minutes for smaller instances to approximately 7–9 hours for the most challenging scenarios, which involve more than 90 product promotion groups (PPGs) over full-year planning horizons and extensive constraint sets. 

To manage computational tractability, the system incorporates advanced callback mechanisms that monitor optimization progress and enable early termination when improvements in the MILP gap stagnate. Target MILP gaps are defined on a market-specific basis, ranging from 1\% in high-value markets where solution precision is critical to 5\% in markets where faster turnaround times are prioritized. Despite the scale and complexity of these problems, the optimization process consistently produces solutions within the prescribed MILP gaps, with results extensively validated by the relevant business units.

\subsection{PricingAI}

PricingAI has been successfully deployed in two major markets, Mexico and the United States, demonstrating its capability to manage complex pricing decisions across diverse and competitive commercial landscapes.

Compared to PromoAI, PricingAI addresses a substantially different class of optimization problems. One of the primary challenges lies in the scale of the product portfolio it manages. The system is designed to optimize prices for over 40 PPGs, encompassing both PepsiCo products and competitor SKUs. 

While the number of constraints involved in PricingAI is generally lower than that in PromoAI, the nature of these constraints is markedly more complex. Specifically, the system must navigate highly nonlinear constraint structures that arise from the logarithmic and exponential components in the volume equations. Furthermore, it incorporates ratio-based functions to model competitive market share, adding another layer of mathematical intricacy.

The decision space also varies by market. In the United States the model operates over continuous price variables (Dollars), whereas in Mexico it deals with integer-valued prices (Pesos). In both cases, psychological pricing thresholds must be respected, introducing additional considerations into the formulation of feasible and optimal price points.

Unlike MILP methods, where the optimization process is deterministic and guarantees convergence to a global optimum given sufficient time and resources, DE is a heuristic, population-based algorithm with no such guarantee. DE relies on stochastic processes for exploration and exploitation, which introduces variability in the results and a risk of convergence to local optima. 

To address this stochastic nature and lack of global optimality guarantees, we employ a strategy of parallel execution consisting of up to 10 independent DE runs. Each run is initialized with distinct algorithmic configurations, including varied mutation factors, crossover probabilities, and population sizes. This ensemble approach enhances the robustness of the optimization process by promoting diverse exploration of the solution space and mitigating the risk of premature convergence to local optima. The optimization tasks are distributed across the Databricks cloud computing infrastructure as described before, allowing concurrent execution of the DE instances and significantly reducing the overall computational time. This design leverages algorithmic diversity and parallelism to improve solution quality and scalability in complex optimization scenarios. 

To ensure reproducibility and deterministic behavior, we incorporate fixed seeds for the DE algorithm. By explicitly setting the seed for the random number generators used in population initialization, mutation, and crossover operations, we guarantee that repeated executions with identical configurations and data yield the same results. This deterministic setup is particularly valuable for benchmarking, debugging, and validating optimization performance, as well as for business coherence, since it eliminates variability due to stochastic components. Combined with our strategy of parallel executions and parameter variation, the use of seeded runs provides a controlled yet comprehensive exploration of the solution space.

The optimization problem takes an average computational time of less than 30 minutes per instance. For an accurate stopping criterion, a strict optimality tolerance is defined. In the differential evolution algorithm, this tolerance parameter governs convergence by monitoring the relative difference between successive population energies. Setting the tolerance to a very low value ensures coherent and reliable solutions while keeping runtimes within the 30-minute limit. Although a maximum runtime is also enforced as a safeguard, the solver typically terminates when the tolerance condition is satisfied, providing a controlled balance between solution quality and computational feasibility.

\subsection{Solution validation}

Before full-scale deployment, PepsiCo conducted both technical and business validation to ensure that the forecasting, elasticity estimation, and optimization components performed as intended and produced recommendations consistent with managerial objectives.

The forecasting models were evaluated using rolling-origin backtesting on historical data. Performance metrics such as MAPE were computed at the SKU–retailer–week level and compared against legacy forecasting benchmarks. For PricingAI, estimated price elasticities were additionally reviewed for economic consistency, ensuring negative own-price effects and reasonable cross-price relationships across products. Elasticity magnitudes were compared against historical pricing cycles and managerial expectations to confirm plausibility. The optimization components of both systems were verified by confirming that all business constraints were satisfied under historical and simulated planning scenarios. Sensitivity analyses were conducted to assess the stability of recommended prices and promotional calendars under moderate variations in demand parameters, cost inputs, and trade spending assumptions.

During this validation phase, planners and revenue management teams reviewed model outputs through structured scenario analyses. Recommended price points and promotional plans were compared with historical decisions to identify systematic differences. When recommendations deviated from expected promotional patterns, the team examined whether the differences were driven by forecast updates (including the treatment of outliers), interactions among constraints, or optimization trade-offs between volume and margin objectives. This iterative review process strengthened user confidence in the system and led to refinements in several business rules embedded in the model. Both solutions were initially deployed in a limited pilot involving selected retailers, with a phased plan to expand to additional markets following successful validation. Throughout the pilot, user feedback and override behavior were monitored to identify potential model adjustments prior to broader rollout.

Sustaining performance over time requires ongoing model maintenance as market conditions, product assortments, and retailer behaviors evolve. The demand models are retrained periodically as new data become available — quarterly for promotional planning and at each pricing cycle for base pricing — following the same validation gates applied at initial deployment before any updated model is released for planning use. Optimization parameters are reviewed each planning cycle to reflect changes in costs, trade spending agreements, retailer requirements, and strategic objectives. Both systems are also designed for extensibility: new product categories and retailers can be incorporated without redesigning the underlying models, lowering the cost of expanding scope over time. In PricingAI, sequential prior updating provides an additional source of stability across refreshes, as posterior estimates from a previous training run inform priors in the next, reducing the risk of large destabilizing shifts in elasticity estimates between planning cycles.

\section{User Experience and Impact}\label{sec:Results}

A critical factor in the success of both systems has been the design of human-in-the-loop workflows that give planners transparency, control, and the ability to incorporate business judgment that the model cannot capture. Rather than presenting optimization outputs as black-box recommendations, both tools are built around an interactive scenario-based interface that supports iterative dialogue between algorithmic recommendations and business expertise.

\subsection{PromoAI: User Workflow}

PromoAI integrates within end-to-end commercial planning workflows through an iterative, retailer-facing process. The workflow begins with planners collecting business rules and constraints directly from the retailer — including budgets, promotional frequency requirements, execution windows, and product exclusivity conditions. These are encoded as model constraints and the optimizer generates multiple promotional calendars under different objective function configurations, for example prioritizing revenue maximization in one scenario and margin preservation in another. This set of optimized calendars is then presented to the retailer as a menu of analytically grounded options, supporting a structured commercial negotiation. 

Once retailer feedback is gathered, planners incorporate the updated requirements — revised constraints, adjusted targets, or new business rules — and re-run the optimization to generate a refined calendar that reflects the agreed commercial framework. This closed-loop interaction between the optimizer and the retailer ensures that the final promotional calendar reflects both analytical rigor and commercial practicality. 

Beyond this iterative optimization workflow, the system also functions as a simulation tool. Once a calendar has been generated, planners can apply manual modifications to individual promotional decisions and immediately observe the predicted impact on all relevant key performance indicators (revenue, margin, volume, and promotional frequency) across the full planning horizon. This capability is particularly valuable for incorporating qualitative considerations that the model cannot anticipate, such as last-minute retailer requests or competitive responses, while maintaining full visibility of the financial consequences of any departure from the optimized plan.

\subsection{PricingAI: User Workflow}

PricingAI mirrors this human-centered design through a dedicated interface that provides planners with full transparency and control over pricing recommendations. The workflow begins with the system generating an optimized price vector across the product portfolio given a set of user-defined constraints: minimum margin thresholds, revenue targets, volume bounds, and price ladder requirements. Planners can then simulate alternative scenarios by modifying these constraints and generating new optimized recommendations in real time. 

A particularly important capability is the adjusted scenario feature, in which users can manually modify individual prices and immediately observe the predicted impact on volume, revenue, and margin across the full portfolio, including cross-price effects on competing SKUs. This functionality has proven especially useful for incorporating qualitative factors not captured by the model, such as retailer relationship considerations or responses to recent competitive moves. The ability to see the full downstream consequences of a manual price change, rather than simply overriding a recommendation, has been central to building organizational trust in the system.

\subsection{Adoption and Business Impact}

The flexibility and transparency of both interfaces have contributed significantly to widespread adoption. PromoAI is now active across all global markets, supporting 16 promotional calendar refreshes per year on a regular planning cadence, while PricingAI has been deployed across US and Mexico markets with further global expansion underway across other markets.

Adoption has translated into strong and measurable operational outcomes. In most markets, about 85\% of PromoAI's optimized promotional recommendations were accepted and executed by the business — a meaningful signal of organizational trust, given that value is only realized when recommendations are actually implemented rather than simply generated. Qualitative feedback from commercial teams has been consistent with these outcomes. A sales director at a major retailer summarized the experience as follows: \textit{``The output is really, really good and allows us to look at and engage with different recommendations. As we build the calendar, this is exactly what we need''}. More broadly, the shift from spreadsheet‑driven planning to optimization‑supported workflows has compressed planning and execution cycles from weeks to minutes for specific pricing and promotional tasks, enabling planners to generate, test, and refine multiple scenarios within a single working session rather than across multiple planning iterations. These efficiency gains have allowed markets to scale usage across retailers and categories without proportional increases in planning effort or headcount, materially lowering the barrier to adoption in new geographies. 

The business impact of both systems has been independently confirmed 
by PepsiCo senior management, who verified that PromoAI and PricingAI 
have delivered measurable improvements in revenue and margin performance, 
alongside more accurate promotional and pricing decisions and reduced 
planning cycle times across all deployed markets. While specific 
financial figures are commercially sensitive and cannot be disclosed, 
industry benchmarks suggest that AI-enabled pricing tools can drive 
2.0--5.0\% improvements in both revenue and profit 
\citep{bcg2021aipricing}.

This transformation extends beyond revenue or technology adoption to encompass organizational culture. Revenue management teams now expect data-driven insights and optimization support as standard capabilities rather than special exceptions. Planning meetings that once revolved around spreadsheet reviews now focus on scenario analysis and strategic discussions enabled by optimization results. What began as a pilot project in a single market has grown into a global capability that influences billions of dollars in revenue decisions, illustrating how patient investment in OR applications can yield transformative organizational results. 

\section{Discussion}\label{sec:Discussion}

The successful deployment of PromoAI and PricingAI across multiple global markets provides valuable insights into the practical implementation of advanced optimization techniques in enterprise CPG environments. Table~\ref{tab:comparison} summarizes the key differences in computational characteristics between the two systems. These differences are not incidental — they reflect the fundamentally different mathematical structures of the two problems. PromoAI's discrete, linearly-constrained scheduling problem is naturally suited to MILP, which provides bounded optimality gaps of 1--5\% but scales to hundreds of thousands of constraints and runtimes of up to 10 hours. PricingAI's nonlinear, non-convex landscape makes exact methods intractable, motivating the use of DE across multiple parallel instances — trading optimality guarantees for tractability, while consistently producing solutions accepted by business stakeholders in practice.

\begin{table}[htbp]
\centering
\caption{Computational characteristics comparison of both projects}
\label{tab:comparison}
\begin{tabular}{lcc}
\hline
\textbf{Characteristic} & \textbf{PromoAI} & \textbf{PricingAI} \\
\hline
Optimization Type & MILP & Nonlinear (DE) \\
Variable Type & Binary & Continuous and Integer \\
Constraint Count & Hundreds of thousands & Dozens \\
Constraint Type & Linear & Nonlinear \\
Typical Runtime & Minutes to 10 hours & Half an hour \\
Parallelization & Single solver instance & Up to 10 parallel runs \\
Markets Deployed & 7+ countries & 2 countries \\
\hline
\end{tabular}
\end{table}

\subsection{The Relative Contributions of Forecasting and Optimization}

A natural question is whether the primary source of business value in these systems lies in the ML forecasting component or in the mathematical optimization component. This section addresses that question directly.

In PromoAI, the case for optimization as the primary value driver is unambiguous. The business problem is fundamentally a large-scale combinatorial scheduling problem: even with perfect demand forecasts, no planner could manually construct promotional calendars that simultaneously satisfy the full set of business constraints — financial targets, calendar rules, exclusivity requirements, and spacing constraints — across hundreds of PPGs and dozens of planning weeks. The feasible solution space contains millions of product-promotion-timing combinations, and the interaction effects between simultaneous promotions make greedy or heuristic manual approaches unreliable. The MILP optimization layer is what makes this tractable: it systematically navigates this space to find schedules that are both financially optimal and operationally valid, with bounded optimality gaps of 1--5\%. The LightGBM forecasting model provides the lift estimates that parameterize the optimization, but without the optimizer these forecasts would support only ad hoc manual planning — valuable, but far from the systematic, constraint-satisfying calendar generation that PromoAI delivers.

In PricingAI, the optimization argument is equally compelling. Base pricing decisions involve setting continuous price changes across a portfolio of 40+ interrelated products, where every individual price movement triggers cross-price effects on all other products in the portfolio simultaneously. A planner adjusting a single product's price by even a fraction of a percent must account for how that change ripples through substitution and complementarity relationships across the full assortment — an interdependency structure that is impossible to reason about manually at portfolio scale. The DE optimization navigates this continuous, non-convex landscape to find the price vector that maximizes the chosen KPI across the full portfolio while simultaneously satisfying price ladder constraints, financial bounds, and all cross-price interaction effects. The Bayesian hierarchical model provides the elasticity matrix that makes this computation possible — a necessary input, but one that alone yields no pricing recommendation. Without the optimization layer, the elasticity model would be a descriptive analytics tool; it is the optimizer that makes it prescriptive and capable of identifying pricing configurations that no manual process could reliably construct.

In both systems, forecasting and optimization are complementary and neither is sufficient alone: forecasting without optimization leaves planners with better information but no mechanism for translating it into optimal decisions at scale, while optimization without reliable forecasts would produce precisely optimal solutions to the wrong problem. The forecasting components are carefully chosen to meet the specific demands of each problem — accuracy, uncertainty quantification, scalability — but they are best understood as enabling inputs to the optimization rather than as ends in themselves. It is the optimization layer that converts predictive analytics into prescriptive value.

\subsection{Implementation Challenges}

Deploying optimization systems at enterprise scale involves challenges that extend well beyond algorithmic design. Both systems are fundamentally dependent on the quality and consistency of historical data, and this proved to be one of the most persistent challenges in practice. CPG retail data varies significantly across markets in terms of granularity, completeness, and reliability. Promotional flagging, price derivation, and SKU matching all require market-specific calibration, and assumptions that work well in one market may fail entirely in another. This heterogeneity means that deployment is not a one-time engineering effort but an ongoing process of market-specific adaptation, particularly as retailer data feeds and product assortments evolve over time.

Several deliberate approximations were necessary to make both systems computationally tractable at production scale. In PromoAI, the piecewise linear approximation of the LightGBM demand forecast converts a non-convex problem into a solvable MILP, reducing solution times from potentially days to hours. The aggregation of individual SKUs into PPGs further reduces dimensionality while preserving decision relevance, reflecting how business users naturally think about products as families rather than individual items. In PricingAI, the use of posterior mean elasticity estimates rather than full uncertainty distributions keeps the optimization tractable. Despite their effectiveness, these approximations introduce real limitations: PromoAI's MILP framework cannot express certain nonlinear business behaviors without workarounds; PricingAI's reliance on metaheuristic optimization means that optimality gaps cannot be formally bounded.

Gaining organizational trust in both the ML forecasts and the optimization outputs was equally challenging, particularly in the early stages of deployment. Planners accustomed to experience-based decision making were often skeptical of recommendations generated by models they could not directly interrogate. This was addressed through the design of transparent, interactive interfaces that allow users to inspect the demand forecasts underlying each recommendation, modify constraints, and simulate the financial consequences of manual overrides. A recurring practical obstacle in PromoAI was that optimized promotional calendars sometimes included recommendations that were commercially sound but operationally infeasible at the retailer level — retailers have their own planning cycles, shelf space constraints, and vendor commitments that are not fully captured by the model's constraint framework. As a result, the iterative workflow in which planners present multiple scenarios to the retailer, gather feedback, and re-optimize is not merely a design preference but a practical necessity: recommendations that cannot be executed have no value, and the system's ability to rapidly regenerate calendars under revised constraints is what makes it commercially viable in practice.

Across both systems, the single most important factor in driving adoption was giving users meaningful control over the optimization process rather than asking them to accept or reject outputs wholesale. The ability to modify recommendations, run scenarios, and see the quantified consequences of departing from the model's suggestion allowed planners to engage with the tools as decision support systems rather than black boxes. This human-in-the-loop design reflects a broader lesson: in enterprise settings, the organizational challenge of gaining trust and changing established workflows is at least as difficult as the technical challenge of building the models themselves, and deployment strategies that underestimate this tend to fail regardless of algorithmic quality.

\section{Conclusion}\label{sec:Conclusion}

This paper has presented two large-scale optimization systems—PromoAI and PricingAI—that illustrate how advanced OR and ML techniques can be applied to revenue management challenges in the CPG industry. Through deployment across multiple markets on several continents, these systems show how sophisticated mathematical optimization can be integrated into enterprise workflows to support large-scale decision-making.

First, we have demonstrated scalable optimization at enterprise scale, showing that MILP and metaheuristic approaches can be used to address the complexity of real-world CPG planning problems. These systems manage hundreds of thousands of constraints and variables while maintaining reasonable solution times, illustrating that academic optimization techniques can be adapted to meet the operational requirements of global business environments.

Second, the systems support a human-centric approach to AI integration. By serving as decision-support tools that augment rather than replace human judgment, they have been adopted by business users and incorporated into existing planning processes. This highlights how advanced analytics can be designed to complement domain expertise rather than fully automate complex commercial decisions.

Third, the use of pragmatic design choices — including piecewise linear 
approximations of demand curves, SKU aggregation into PPGs, and 
JSON-based modular constraint configuration — illustrates how theoretical 
rigor can be balanced with computational efficiency and deployment 
scalability. These techniques enable otherwise intractable problems to 
be solved at scale while accommodating diverse business rules across 
heterogeneous markets, offering practical guidance for practitioners 
facing similar implementation challenges.

In sum, PromoAI and PricingAI illustrate that the path from theory to real-world application is achievable through a combination of technical soundness, business alignment, and thoughtful system design. These systems reinforce the role of OR as a discipline capable of supporting scalable, analytics-driven decision-making and augmenting human judgment in complex business settings.


%


\section*{Appendix A. Mathematical Formulations for PromoAI}

This appendix provides the mathematical details underlying the PromoAI
optimization engine. The main text presents the concepts in natural language
in accordance with IJAA guidelines; the complete mathematical structure is
documented here for transparency. We first declare sets, parameters, and
decision variables; we then state the full optimization model; finally, we
describe each component in turn.

\subsection*{A.1 Forecasting Model}

PromoAI uses a LightGBM regression model \citep{ke2017lightgbm} to estimate
promotional demand. The model predicts weekly unit sales for each Product
Promotion Group (PPG) as a nonlinear function of pricing, calendar variables,
promotional mechanics, and product attributes:

\begin{equation}
\hat{d}_{p,w,r} = f\bigl(\mathrm{price},\; \mathrm{week},\; \mathrm{mechanic},\; \mathrm{PPG},\; \ldots\bigr)
\label{eq:forecast}
\end{equation}

\noindent where $f(\cdot)$ is the function learned through gradient-boosted
decision trees. These forecasts feed the piecewise linear demand approximation
embedded in the optimization model below.

Table \ref{tab:lgbm_hyperparams} summarizes the typical hyperparameter ranges explored and the values selected for representative markets.

\begin{table}[h]
\centering
\caption{Typical LightGBM hyperparameter search ranges and selected values for PromoAI forecasting models across markets. The objective function is MAPE, which is robust to outliers in sales data. Early stopping (patience = 50 rounds) on a held-out validation set determines the final number of trees. Hyperparameters are tuned via randomized cross-validation on a temporal hold-out split.}
\label{tab:lgbm_hyperparams}
\begin{tabular}{llll}
\hline
\textbf{Hyperparameter} & \textbf{Search Range} & \textbf{Selected Value} & \textbf{Role} \\
\hline
\texttt{num\_leaves}              & 50--500                        & 269    & Controls model capacity \\
\texttt{learning\_rate}           & Loguniform $[e^{-4},\, e^{2}]$ & 0.0185 & Step size shrinkage \\
\texttt{n\_estimators}            & 300--600                       & 302    & Number of trees \\
\texttt{min\_split\_gain}         & 0.001--0.01                    & 0.00456 & Min gain to split a node \\
\texttt{max\_depth}               & 100--250                       & 184    & Tree depth limit \\
\texttt{reg\_alpha}               & 0--10                          & 6.04   & L1 regularization \\
\texttt{reg\_lambda}              & 0--10                          & 4.95   & L2 regularization \\
\hline
\end{tabular}
\end{table}

\subsection*{A.2 Sets and Indices}

Table~\ref{tab:sets} lists the sets and index conventions used throughout
the formulation.

\begin{table}[!htbp]
\centering
\caption{Sets used in the PromoAI formulation.}
\label{tab:sets}
\begin{tabular}{cl}
\toprule
\textbf{Set} & \textbf{Description} \\
\midrule
$\mathcal{P}$ & Set of all Product Promotion Groups (PPGs) \\
$\widetilde{\mathcal{P}} \subseteq \mathcal{P}$ & Subset of PepsiCo-owned PPGs \\
$\mathcal{W}$ & Set of weeks in the planning horizon \\
$\mathcal{R}_p$ & Set of promotional options available for PPG $p \in \mathcal{P}$ \\
$\widetilde{\mathcal{R}}_p \subseteq \mathcal{R}_p$ & Subset of non-base (active) promotions for PPG $p$ \\
$\mathcal{S}_p \subseteq \mathcal{P}$ & Set of PPGs in the same size/format segment as $p$ (excluding $p$) \\
$\mathcal{L}$ & Set of locked competitor promotion assignments $(p,w,r)$ \\
\bottomrule
\end{tabular}
\end{table}

\subsection*{A.3 Parameters}

Table~\ref{tab:params} summarizes the parameters of the formulation,
together with their units. All parameters are retrieved from the scenario
configuration and the input data.

\begin{table}[!htbp]
\centering
\caption{Parameters of the PromoAI formulation. All parameters are
         retrieved from the scenario configuration and input data.}
\label{tab:params}
\begin{tabular}{cll}
\toprule
\textbf{Symbol} & \textbf{Description} & \textbf{Units} \\
\midrule
$\delta_{p,r}$
  & Discount fraction for PPG $p$, promotion $r$
  & $(-)$ \\
$\rho^{\mathrm{P}}_{p,w,r}$
  & PepsiCo unit revenue
  & (currency/unit) \\
$\mu^{\mathrm{P}}_{p,w,r}$
  & PepsiCo unit margin
  & (currency/unit) \\
$\rho^{\mathrm{R}}_{p,w,r}$
  & Retailer unit revenue
  & (currency/unit) \\
$\mu^{\mathrm{R}}_{p,w,r}$
  & Retailer unit margin
  & (currency/unit) \\
$\omega_1$
  & Weight on PepsiCo in objective
  & $(-)$ \\
$\omega_2$
  & Weight on retailer in objective
  & $(-)$ \\
$\omega_3$
  & Weight on PepsiCo margin vs.\ sales
  & $(-)$ \\
$\omega_4$
  & Weight on retailer margin vs.\ sales
  & $(-)$ \\
$\underline{\gamma}^{\mathrm{P}}$
  & Min.\ PepsiCo margin-to-sales ratio
  & $(-)$ \\
$\underline{\gamma}^{\mathrm{R}}$
  & Min.\ retailer margin-to-sales ratio
  & $(-)$ \\
$\overline{N}$
  & Max.\ total active promotions
  & (promotions) \\
$\overline{N}^{\mathrm{P}}$
  & Max.\ PepsiCo-only promotions
  & (promotions) \\
$\mathbf{b}_{p,w,r}$
  & Breakpoint vector for PWL demand approximation
  & $(-)$ \\
$\mathbf{d}_{p,w,r}$
  & Demand values at each PWL breakpoint
  & (units) \\
$\underline{m}$
  & Min.\ PepsiCo market share
  & $(-)$ \\
$M$
  & Sufficiently large constant (big-$M$)
  & $(-)$ \\
\bottomrule
\end{tabular}
\end{table}

\subsection*{A.4 Decision and Auxiliary Variables}

Table~\ref{tab:variables} defines the primary decision variable and the
auxiliary variables used in the formulation. For readability, we also define aggregated financial auxiliary variables in Table~\ref{tab:finauxvars}.

\begin{table}[!htbp]
\centering
\caption{Decision and auxiliary variables.}
\label{tab:variables}
\begin{tabular}{clll}
\toprule
\textbf{Symbol} & \textbf{Domain} & \textbf{Description} & \textbf{Units} \\
\midrule
\multicolumn{4}{l}{\textit{Primary decision variable}} \\
$x_{p,w,r}$ & $\{0,1\}$
  & 1 if promo $r$ selected for PPG $p$, week $w$
  & $(-)$ \\
\addlinespace
\multicolumn{4}{l}{\textit{Auxiliary variables}} \\
$\alpha_{p,w}$ & $\mathbb{R}_{\geq 0}$
  & Discount fraction, PPG $p$, week $w$
  & $(-)$ \\
$\beta_{p,w}$ & $\mathbb{R}_{\geq 0}$
  & Discount pressure, PPG $p$, week $w$
  & $(-)$ \\
$\hat{q}_{p,w,r}$ & $\mathbb{R}_{\geq 0}$
  & Conditional demand (if promo active)
  & (units) \\
$q_{p,w,r}$ & $\mathbb{R}_{\geq 0}$
  & Realized demand
  & (units) \\
\bottomrule
\end{tabular}
\end{table}

\begin{table}[!htbp]
\centering
\caption{Aggregated financial auxiliary variables.}
\label{tab:finauxvars}
\begin{tabular}{cll}
\toprule
\textbf{Symbol} & \textbf{Description} & \textbf{Units} \\
\midrule
$S^{\mathrm{P}}_{p,w}$, $\Pi^{\mathrm{P}}_{p,w}$
  & PepsiCo sales and margin (per PPG-week)
  & (currency) \\
$S^{\mathrm{R}}_{p,w}$, $\Pi^{\mathrm{R}}_{p,w}$
  & Retailer sales and margin (per PPG-week)
  & (currency) \\
$S^{\mathrm{P}}$, $\Pi^{\mathrm{P}}$
  & Total PepsiCo sales and margin
  & (currency) \\
$S^{\mathrm{R}}$, $\Pi^{\mathrm{R}}$
  & Total retailer sales and margin
  & (currency) \\
\bottomrule
\end{tabular}
\end{table}

\subsection*{A.5 Optimization Model}

The complete PromoAI optimization model, using the sets, parameters, and
variables defined above, is stated as follows.

\medskip
\noindent\textit{Objective function.}

\begin{equation}
\begin{aligned}
\max \quad & \omega_1 \bigl[(1 - \omega_3) \, S^{\mathrm{P}} + \omega_3 \, \Pi^{\mathrm{P}} \bigr] \\
& +\; \omega_2 \bigl[(1 - \omega_4) \, S^{\mathrm{R}} + \omega_4 \, \Pi^{\mathrm{R}} \bigr]
\end{aligned}
\label{eq:objective}
\end{equation}

\medskip
\noindent\textit{Subject to:}

\medskip
\noindent\textit{Promotion selection.} For all $p \in \mathcal{P}$,
$w \in \mathcal{W}$:

\begin{equation}
\sum_{r \in \mathcal{R}_p} x_{p,w,r} = 1
\label{eq:exclusivity}
\end{equation}

\begin{equation}
\alpha_{p,w} = \sum_{r \in \mathcal{R}_p} \delta_{p,r} \cdot x_{p,w,r}
\label{eq:discount_frac}
\end{equation}

\medskip
\noindent\textit{Demand computation.} For all $p \in \mathcal{P}$,
$w \in \mathcal{W}$, $r \in \mathcal{R}_p$:

\begin{equation}
\beta_{p,w} = \frac{1}{|\mathcal{S}_p|} \sum_{p' \in \mathcal{S}_p} \alpha_{p',w}
\label{eq:discount_pressure}
\end{equation}

\begin{equation}
\hat{q}_{p,w,r} = \mathrm{PWL}\bigl(\beta_{p,w}\,;\; \mathbf{b}_{p,w,r},\; \mathbf{d}_{p,w,r}\bigr)
\label{eq:pwl}
\end{equation}%
\vspace*{-30pt}%
\begin{subequations}
\label{eq:bigm}
\begin{align}
\hat{q}_{p,w,r} - q_{p,w,r} &\leq M \cdot (1 - x_{p,w,r})
\label{eq:bigm_a} \\
q_{p,w,r} &\leq M \cdot x_{p,w,r}
\label{eq:bigm_b} \\
q_{p,w,r} &\leq \hat{q}_{p,w,r}
\label{eq:bigm_c}
\end{align}
\end{subequations}

\newpage
\noindent\textit{Financial aggregation.} For all $p \in \mathcal{P}$,
$w \in \mathcal{W}$:
\vspace*{-20pt}
\begin{subequations}
\label{eq:financial_pw}
\begin{align}
S^{\mathrm{P}}_{p,w} &= \sum_{r \in \mathcal{R}_p} q_{p,w,r} \cdot \rho^{\mathrm{P}}_{p,w,r}
\label{eq:pep_sales_pw} \\
\Pi^{\mathrm{P}}_{p,w} &= \sum_{r \in \mathcal{R}_p} q_{p,w,r} \cdot \mu^{\mathrm{P}}_{p,w,r}
\label{eq:pep_margin_pw} \\
S^{\mathrm{R}}_{p,w} &= \sum_{r \in \mathcal{R}_p} q_{p,w,r} \cdot \rho^{\mathrm{R}}_{p,w,r}
\label{eq:ret_sales_pw} \\
\Pi^{\mathrm{R}}_{p,w} &= \sum_{r \in \mathcal{R}_p} q_{p,w,r} \cdot \mu^{\mathrm{R}}_{p,w,r}
\label{eq:ret_margin_pw}
\end{align}
\end{subequations}%
\vspace*{-50pt}%
\begin{subequations}
\label{eq:financial_totals}
\begin{align}
S^{\mathrm{P}} &= \sum_{p \in \mathcal{P}} \sum_{w \in \mathcal{W}} S^{\mathrm{P}}_{p,w}
\label{eq:pep_sales_total} \\
\Pi^{\mathrm{P}} &= \sum_{p \in \mathcal{P}} \sum_{w \in \mathcal{W}} \Pi^{\mathrm{P}}_{p,w}
\label{eq:pep_margin_total} \\
S^{\mathrm{R}} &= \sum_{p \in \mathcal{P}} \sum_{w \in \mathcal{W}} S^{\mathrm{R}}_{p,w}
\label{eq:ret_sales_total} \\
\Pi^{\mathrm{R}} &= \sum_{p \in \mathcal{P}} \sum_{w \in \mathcal{W}} \Pi^{\mathrm{R}}_{p,w}
\label{eq:ret_margin_total}
\end{align}
\end{subequations}

\medskip
\noindent\textit{Business constraints.}
\vspace*{-20pt}
\begin{subequations}
\label{eq:margin}
\begin{align}
\Pi^{\mathrm{P}} &\geq \underline{\gamma}^{\mathrm{P}} \cdot S^{\mathrm{P}}
\label{eq:margin_pep} \\
\Pi^{\mathrm{R}} &\geq \underline{\gamma}^{\mathrm{R}} \cdot S^{\mathrm{R}}
\label{eq:margin_ret}
\end{align}
\end{subequations}%
\vspace*{-50pt}%
\begin{subequations}
\label{eq:promo_caps}
\begin{align}
\sum_{p \in \mathcal{P}} \sum_{w \in \mathcal{W}} \sum_{r \in \widetilde{\mathcal{R}}_p} x_{p,w,r} &\leq \overline{N}
\label{eq:total_promo_cap} \\
\sum_{p \in \widetilde{\mathcal{P}}} \sum_{w \in \mathcal{W}} \sum_{r \in \widetilde{\mathcal{R}}_p} x_{p,w,r} &\leq \overline{N}^{\mathrm{P}}
\label{eq:pep_promo_cap}
\end{align}
\end{subequations}%
\vspace*{-30pt}
\begin{equation}
\sum_{p \in \widetilde{\mathcal{P}}} \sum_{w \in \mathcal{W}} S^{\mathrm{R}}_{p,w}
\;\geq\;
\underline{m} \cdot S^{\mathrm{R}}
\label{eq:market_share}
\end{equation}

\begin{equation}
x_{p,w,r} = 1 \qquad \forall\, (p,w,r) \in \mathcal{L}
\label{eq:competitor_lock}
\end{equation}

\newpage
\noindent\textit{Variable domains.}

\begin{equation}
\begin{aligned}
& x_{p,w,r} \in \{0,1\} \\
& \alpha_{p,w},\; \beta_{p,w},\; \hat{q}_{p,w,r},\; q_{p,w,r} \in \mathbb{R}_{\geq 0} \\
& \forall\, p \in \mathcal{P},\; w \in \mathcal{W},\; r \in \mathcal{R}_p
\end{aligned}
\label{eq:promo_domains}
\end{equation}

\subsection*{A.6 Description of Model Components}

This subsection describes each component of the model stated in
Section~A.5.

\subsubsection*{A.6.1 Objective Function}

The objective~\eqref{eq:objective} maximizes a user-configurable weighted
combination of PepsiCo and retailer financial performance, subject to
$\omega_1 + \omega_2 = 1$. The four weights are specified by the user
through the system interface. For example, setting $\omega_1 = 1$ and
$\omega_3 = 0$ reduces the objective to maximizing PepsiCo sales only.

\subsubsection*{A.6.2 Promotion Exclusivity}

Constraint~\eqref{eq:exclusivity} requires each PPG to be assigned exactly
one promotional option (including the no-promotion baseline) in every week.

\subsubsection*{A.6.3 Discount Fraction}

Constraint~\eqref{eq:discount_frac} computes the effective discount for each
PPG-week from the selected promotion.

\subsubsection*{A.6.4 Competitive Discount Pressure}

Constraint~\eqref{eq:discount_pressure} captures the average promotional
intensity of competing products within the same size or format segment.
This auxiliary variable feeds the piecewise linear demand approximation,
encoding a measure of competitive cannibalization into the demand forecast.

\subsubsection*{A.6.5 Piecewise Linear Demand Approximation}

The machine learning forecast $f(\cdot)$ produces nonlinear demand responses,
particularly when promotions interact with competitor discount levels. To
maintain a tractable MILP formulation, these nonlinear response surfaces
are approximated using piecewise linear (PWL) functions. Depending on the
market, two to five breakpoints are selected to balance accuracy with
computational efficiency.

In~\eqref{eq:pwl}, the conditional demand $\hat{q}_{p,w,r}$ is defined
through a PWL function of the discount pressure, where $\mathbf{b}_{p,w,r}$
and $\mathbf{d}_{p,w,r}$ are vectors of breakpoints (discount pressure
values) and corresponding demand values, respectively, derived from the
ML forecast during scenario generation.

\subsubsection*{A.6.6 Demand Activation via Big-$M$ Linearization}

The realized demand $q_{p,w,r}$ must equal the conditional demand
$\hat{q}_{p,w,r}$ when promotion $r$ is selected (i.e.,
$x_{p,w,r} = 1$) and must be zero otherwise. The big-$M$
constraints~\eqref{eq:bigm} enforce this implication.
Constraint~\eqref{eq:bigm_a} ensures that when $x_{p,w,r} = 1$, the
realized demand is at least as large as the conditional demand.
Constraint~\eqref{eq:bigm_b} forces the realized demand to zero when
the promotion is not selected. Constraint~\eqref{eq:bigm_c} provides
an upper bound. Together with the maximization objective, which pushes
demand upward, these three constraints enforce
$q_{p,w,r} = \hat{q}_{p,w,r} \cdot x_{p,w,r}$ at optimality.

\subsubsection*{A.6.7 Financial Intermediate Calculations}

Equations~\eqref{eq:financial_pw} compute per-PPG, per-week financial
metrics from the realized demand and per-unit financial parameters.
Equations~\eqref{eq:financial_totals} aggregate these over all PPGs and
weeks to obtain the totals used in the objective function.

\subsubsection*{A.6.8 Margin Conservation}

Constraints~\eqref{eq:margin} enforce minimum margin-to-sales ratios
for both PepsiCo and the retailer. The ratio constraints are linearized
by multiplying through by the denominator, which is strictly positive in
any feasible solution.

\subsubsection*{A.6.9 Promotion Frequency Limits}

Constraint~\eqref{eq:total_promo_cap} bounds the total number of active
(non-base) promotions across all PPGs and weeks, while
\eqref{eq:pep_promo_cap} applies the analogous limit to PepsiCo-owned
PPGs only.

\subsubsection*{A.6.10 Market Share}

Constraint~\eqref{eq:market_share} requires PepsiCo's share of total
retailer sales to meet a minimum threshold.

\subsubsection*{A.6.11 Competitor Lock Slots}

Constraint~\eqref{eq:competitor_lock} fixes promotional slots designated
to competitor products in the benchmark calendar.

\subsubsection*{A.6.12 Additional Configurable Constraints}

The constraints presented above represent the core formulation. In
practice, PromoAI supports a broader library of modular constraint
templates that are activated dynamically based on the market
configuration. These include, among others:

\begin{itemize}
\item \textit{Ad-block linking}: constraints that enforce joint promotion
of related products (e.g., requiring complementary pack sizes to be
promoted together in the same week).
\item \textit{Minimum and maximum promotion duration}: constraints on
the number of consecutive weeks a promotion can run.
\item \textit{Weekly promotion caps by mechanic type}: limits on the
number of promotions of a specific type (e.g., price reductions) active
in any given week.
\item \textit{Seasonal and holiday alignment}: constraints that lock or
exclude promotions during specific calendar periods.
\item \textit{Minimum spacing between promotions}: constraints that
enforce a gap of at least $k$ weeks between successive promotions for a
given PPG, preventing consumer pantry-loading effects.
\end{itemize}

\noindent These constraint templates are encoded in JSON-based
configuration files (see Appendix~B) and are instantiated with
market-specific parameters at runtime, allowing business users to
compose complex scenarios without modifying the underlying optimization
code.


\section*{Appendix B. JSON example illustrating two constraints}

Both projects, PromoAI and PricingAI, dynamically constructs the optimization
model by decoding JSON specifications from the web interface. An example of
that for PromoAI can be seen in Code~\ref{json-example}. Abstract building
blocks like MinMaxPromo can be instantiated with different parameters to
create specific business rules across all three constraint categories,
enabling business users to configure complex scenarios without programming
knowledge.

\begin{codefloat}
\begin{lstlisting}[
    frame=single,
    numbers=left,
    xleftmargin=10pt,
    basicstyle=\footnotesize\ttfamily,
    tabsize=4
]
"WeeklyMaxPromo": [
    {
      "data": {
        "mechanic_type": [
          "Price"
        ]
      },
      "week_max_promo": 4
    }
],
"MinMaxPromoDuration": {
    "data":  [
         {
          "pepsico_flag": true
         }
    ],
    "max_promo_block": 4,
    "min_promo_block": 2
},
\end{lstlisting}
\caption{
\textbf{JSON example illustrating two constraints: \texttt{WeeklyMaxPromo} and \texttt{MinMaxPromoDuration}.}
The \texttt{WeeklyMaxPromo} constraint limits the number of weekly promotions of type \texttt{"Price"} to a maximum of 4.
The \texttt{MinMaxPromoDuration} constraint enforces that all promotions for PepsiCo products have a duration between 2 and 4 weeks.
The \texttt{data} field within the JSON supports flexible filtering based on columns from the input dataset.
In this example, filters are applied to the \texttt{mechanic\_type} and \texttt{pepsico\_flag} columns.
This modular structure allows the same codebase to support a wide variety of constraint types through simple configuration changes.
}
\label{json-example}
\end{codefloat}


\section*{Appendix C. Mathematical Formulation for PricingAI}

This appendix provides the mathematical details underlying the PricingAI
optimization engine. The main text presents the concepts in natural
language in accordance with IJAA guidelines; the complete mathematical
structure is documented here for transparency. We first declare sets,
parameters, and decision variables; we then state the full optimization
model as a chain of definitional equations terminating in the objective
and constraints; finally, we describe each component in turn.

\subsection*{C.1 Demand Model}

PricingAI uses a Bayesian hierarchical model to estimate own- and
cross-price elasticities across the product portfolio. The model is
implemented in STAN \citep{carpenter2017stan} using variational
inference (ADVI) \citep{kucukelbir2017advi}. The posterior mean
elasticity estimates are assembled into a cross-price elasticity matrix
$\mathbf{E} \in \mathbb{R}^{n \times n}$, where entry $E_{ij}$ captures
the percentage change in demand for product~$i$ associated with a
one-percent change in the price of product~$j$. These point estimates
feed directly into the volume function defined in the optimization
model below.

\subsection*{C.2 Sets and Indices}

Table~\ref{tab:pricing_sets} lists the sets and index conventions used
throughout the formulation.

\begin{table}[!htbp]
\centering
\caption{Sets used in the PricingAI formulation.}
\label{tab:pricing_sets}
\begin{tabular}{cl}
\toprule
\textbf{Set} & \textbf{Description} \\
\midrule
$\mathcal{A} = \{1,\dots,n\}$
  & Full set of products (PepsiCo and competitors) \\
$\mathcal{P} \subseteq \mathcal{A}$
  & Subset of PepsiCo products \\
$\mathcal{C} \subseteq \mathcal{A}$
  & Subset of competitor products with price-follow behavior \\
$\mathcal{G}$
  & Set of pricing groups at the optimization granularity level \\
$\mathcal{L}_k$
  & Set of products belonging to pricing line $k$ \\
$\mathcal{B}_b$
  & Set of products belonging to brand unit $b$ \\
$\mathcal{T}_t$
  & Set of products in pricing tier $t$ \\
$\mathcal{U}(t)$
  & Set of upper-bound tiers associated with tier $t$ \\
$\mathcal{K}_i$
  & Set of psychological price point (PPP) thresholds for product $i$ \\
\bottomrule
\end{tabular}
\end{table}

\noindent The optimizer operates at the \emph{pricing group} level. A
mapping $g(i): \mathcal{A} \to \mathcal{G}$ associates each product~$i$
with its pricing group, so that products sharing a group receive the
same price decision.

\subsection*{C.3 Parameters}

Table~\ref{tab:pricing_params} summarizes the parameters of the
formulation, grouped by role. All parameters are retrieved from the
scenario configuration and the input data.

\begin{table}[!htbp]
\centering
\caption{Parameters of the PricingAI formulation. All parameters are
         retrieved from the scenario configuration and input data.}
\label{tab:pricing_params}
\begin{tabular}{cll}
\toprule
\textbf{Symbol} & \textbf{Description} & \textbf{Units} \\
\midrule
\multicolumn{3}{l}{\textit{Product parameters}} \\
$p_i^{0}$
  & Baseline shelf price of product $i$
  & (currency/unit) \\
$V_i^{0}$
  & Baseline sales volume of product $i$
  & (units) \\
$s_i$
  & Pack size of product $i$
  & (weight or volume) \\
$c_i^{\mathrm{SIP}}$
  & Sell-in price (manufacturer to retailer)
  & (currency/unit) \\
$d_i$
  & Sell-in discount fraction
  & $(-)$ \\
$c_i^{\mathrm{COGS}}$
  & Cost of goods sold per unit
  & (currency/unit) \\
$c_i^{\mathrm{Dist}}$
  & Distribution cost per unit
  & (currency/unit) \\
$\phi$
  & Sell-in price pass-through factor
  & $(-)$ \\
\addlinespace
\multicolumn{3}{l}{\textit{Elasticity parameters}} \\
$\mathbf{E} \in \mathbb{R}^{n \times n}$
  & Cross-price elasticity matrix
  & $(-)$ \\
$\tau_{i,k}$
  & Psychological price point (PPP) threshold $k$ for product $i$
  & (currency/unit) \\
$\psi_{i,k}$
  & PPP coefficient for threshold $k$ of product $i$
  & $(-)$ \\
\addlinespace
\multicolumn{3}{l}{\textit{Competitor follow parameters}} \\
$\rho_i$
  & Competitor follow ratio for product $i \in \mathcal{C}$
  & $(-)$ \\
$h(i)$
  & PepsiCo ``hero'' product that competitor $i$ follows
  & $(-) $ \\
\addlinespace
\multicolumn{3}{l}{\textit{Constraint parameters}} \\
$\underline{\alpha}^{K}$, $\overline{\alpha}^{K}$
  & Min/max ratio relative to baseline for KPI $K$
  & $(-)$ \\
$\underline{\theta}$, $\overline{\theta}$
  & Min/max average price change (as percentage)
  & $(\%)$ \\
$\gamma$
  & Ladder constraint tolerance factor
  & $(-)$ \\
\bottomrule
\end{tabular}
\end{table}

\subsection*{C.4 Decision Variables}

Table~\ref{tab:pricing_vars} defines the single decision variable of the
PricingAI formulation. The interpretation of $x_g$ depends on the market
configuration: in \textit{percentage mode} (e.g., the U.S.\ market),
$x_g \in \mathbb{R}$ is a multiplicative factor on the baseline price;
in \textit{absolute mode} (e.g., the Mexico market), $x_g \in \mathbb{Z}$
is the new shelf price expressed directly in currency units. Both cases
are formalized in the optimization model in Section~C.5.

\begin{table}[!htbp]
\centering
\caption{Decision variables of the PricingAI formulation.}
\label{tab:pricing_vars}
\begin{tabular}{clll}
\toprule
\textbf{Symbol} & \textbf{Domain} & \textbf{Description} & \textbf{Units} \\
\midrule
$x_g$ & $\mathbb{R}$ or $\mathbb{Z}$
  & Price decision for pricing group $g \in \mathcal{G}$
  & market-dependent \\
\bottomrule
\end{tabular}
\end{table}

\subsection*{C.5 Optimization Model}

The complete PricingAI optimization model is stated below. Auxiliary
quantities (final prices, predicted volumes, financial KPIs) are
introduced as definitional equations and feed into the objective and
constraint expressions.

\medskip
\noindent\textit{Objective function.}

\begin{equation}
\max_{\mathbf{x}} \quad K_{\mathcal{P}}(\mathbf{x})
\label{eq:pricing_obj}
\end{equation}

\noindent where
$K \in \{R,\; \Pi,\; \mathrm{Vol},\; \mathrm{Margin},\; \mathrm{MS},\; \mathrm{NRW}\}$
is selected by the user through the system interface from the KPIs defined
in~\eqref{eq:revenue}--\eqref{eq:nrp}.

\medskip
\noindent\textit{Subject to:}

\medskip
\noindent\textit{Variable bounds.}

\begin{equation}
\underline{x}_g \leq x_g \leq \overline{x}_g \qquad \forall\, g \in \mathcal{G}
\label{eq:var_bounds}
\end{equation}

\medskip
\noindent\textit{Price computation.}
For all $i \in \mathcal{P}$, in percentage mode:

\begin{equation}
p_i^{\mathrm{new}} = x_{g(i)} \cdot p_i^{0}
\label{eq:price_pct}
\end{equation}

\noindent or, in absolute mode:

\begin{equation}
p_i^{\mathrm{new}} = x_{g(i)}
\label{eq:price_abs}
\end{equation}

\noindent For all $i \in \mathcal{P}$:

\begin{equation}
\Delta_i = \frac{p_i^{\mathrm{new}} - p_i^{0}}{p_i^{0}}
\label{eq:pep_pctchange}
\end{equation}

\noindent For all $i \in \mathcal{C}$:

\begin{equation}
\Delta_i = \rho_i \cdot \Delta_{h(i)}
\label{eq:comp_follow}
\end{equation}

\noindent For all $i \in \mathcal{A}$:

\begin{equation}
p_i^{\mathrm{f}} = (1 + \Delta_i) \cdot p_i^{0}
\label{eq:final_price}
\end{equation}

\medskip
\noindent\textit{Volume function.}
For all $i \in \mathcal{A}$:

\begin{subequations}
\label{eq:volume_components}
\begin{align}
\mathcal{E}_i(\mathbf{p}) &= \sum_{j \in \mathcal{A}}
  E_{ij} \bigl[\ln(p_j) - \ln(p_j^{0})\bigr]
\label{eq:elast_component} \\
\mathcal{D}_i(\mathbf{p}) &= \mathrm{PPP}_i(\mathbf{p})
  - \mathrm{PPP}_i(\mathbf{p}^{0})
\label{eq:ppp_component}
\end{align}
\end{subequations}%
\begin{equation}
\mathrm{PPP}_i(\mathbf{p}) = \!\!\sum_{k \in \mathcal{K}_i:\; p_i \unrhd \tau_{i,k}} \!\!\psi_{i,k}
\label{eq:ppp}
\end{equation}

\begin{equation}
V_i(\mathbf{p}^{\mathrm{f}}) = V_i^{0} \cdot
  \exp\!\bigl(\mathcal{E}_i(\mathbf{p}^{\mathrm{f}})
  + \mathcal{D}_i(\mathbf{p}^{\mathrm{f}})\bigr)
\label{eq:volume}
\end{equation}

\medskip
\noindent\textit{Financial KPIs.}
For all $i \in \mathcal{A}$:

\begin{equation}
U_i = \frac{V_i(\mathbf{p}^{\mathrm{f}})}{s_i}
\label{eq:units}
\end{equation}

\begin{equation}
c_i^{\mathrm{SIP,new}} = c_i^{\mathrm{SIP}} \cdot (1 + \Delta_i \cdot \phi)
\label{eq:sip_update}
\end{equation}

\noindent For a product subset $\mathcal{S} \subseteq \mathcal{A}$
(typically $\mathcal{S} = \mathcal{P}$):

\begin{equation}
R_{\mathcal{S}} = \sum_{i \in \mathcal{S}} U_i \cdot c_i^{\mathrm{SIP,new}} \cdot (1 - d_i)
\label{eq:revenue}
\end{equation}

\begin{equation}
\Pi_{\mathcal{S}} = \sum_{i \in \mathcal{S}} \bigl[
  R_i - (c_i^{\mathrm{COGS}} + c_i^{\mathrm{Dist}}) \cdot U_i
\bigr]
\label{eq:profit}
\end{equation}

\noindent where $R_i = U_i \cdot c_i^{\mathrm{SIP,new}} \cdot (1 - d_i)$
is the per-product revenue from~\eqref{eq:revenue}.

\begin{equation}
\mathrm{Margin}_{\mathcal{S}} = \frac{\Pi_{\mathcal{S}}}{R_{\mathcal{S}}}
\label{eq:margin_pricing}
\end{equation}

\begin{equation}
\mathrm{Vol}_{\mathcal{S}} = \sum_{i \in \mathcal{S}} V_i(\mathbf{p}^{\mathrm{f}})
\label{eq:total_volume}
\end{equation}

\begin{equation}
\mathrm{MS}_{\mathcal{S}} = \frac{
  \sum_{i \in \mathcal{S}} U_i \cdot p_i^{\mathrm{f}}
}{
  \sum_{i \in \mathcal{A}} U_i \cdot p_i^{\mathrm{f}}
}
\label{eq:market_share_pricing}
\end{equation}

\begin{equation}
\mathrm{NRW}_{\mathcal{S}} = \frac{R_{\mathcal{S}}}{\mathrm{Vol}_{\mathcal{S}}}
\label{eq:nrp}
\end{equation}

\medskip
\noindent\textit{Financial KPI bounds.}
For each KPI function $K(\mathbf{x})$ with baseline value $K^{0}$:

\begin{equation}
\underline{\alpha}^{K} \cdot K^{0} \;\leq\; K(\mathbf{x}) \;\leq\; \overline{\alpha}^{K} \cdot K^{0}
\label{eq:kpi_bounds}
\end{equation}

\medskip
\noindent\textit{Structural pricing rules.}
For each brand unit $\mathcal{B}_b$, with products sorted by pack size
$s_i$ in descending order, for consecutive products $(i, i+1)$:

\begin{equation}
\frac{p_i^{\mathrm{new}}}{s_i} \;\geq\; \gamma \cdot \frac{p_{i+1}^{\mathrm{new}}}{s_{i+1}}
\label{eq:ladder}
\end{equation}

\noindent For each tier $t$ and its associated upper-bound tiers
$\mathcal{U}(t)$:

\begin{equation}
p_u^{\mathrm{new}} \;\geq\; p_t^{\mathrm{new}} \qquad \forall\, u \in \mathcal{U}(t)
\label{eq:tier}
\end{equation}

\noindent For all products $i, j \in \mathcal{L}_k$:

\begin{equation}
|p_i^{\mathrm{new}} - p_i^{0}| = |p_j^{\mathrm{new}} - p_j^{0}|
\label{eq:pricing_line}
\end{equation}

\medskip
\noindent\textit{Portfolio-level price constraints.}

\begin{equation}
\underline{\theta} \;\leq\;
\frac{\sum_{i \in \mathcal{P}} (p_i^{\mathrm{new}} - p_i^{0})}
     {\sum_{i \in \mathcal{P}} p_i^{0}}
\times 100
\;\leq\; \overline{\theta}
\label{eq:avg_price}
\end{equation}%
\vspace*{-\baselineskip}%
\begin{subequations}
\label{eq:wtd_price}
\begin{align}
\bar{p}^{0} &= \frac{\sum_{i \in \mathcal{P}} p_i^{0} \cdot V_i^{0}}
                     {\sum_{i \in \mathcal{P}} V_i^{0}}
\label{eq:wtd_base} \\[6pt]
\bar{p}^{\mathrm{new}} &= \frac{
  \sum_{i \in \mathcal{P}} p_i^{\mathrm{new}} \cdot V_i(\mathbf{p}^{\mathrm{f}})
}{
  \sum_{i \in \mathcal{P}} V_i(\mathbf{p}^{\mathrm{f}})
}
\label{eq:wtd_new}
\end{align}
\end{subequations}%
\begin{equation}
\underline{\theta} \;\leq\;
\frac{\bar{p}^{\mathrm{new}} - \bar{p}^{0}}{\bar{p}^{0}} \times 100
\;\leq\; \overline{\theta}
\label{eq:wtd_price_bound}
\end{equation}

\subsection*{C.6 Description of Model Components}

This subsection describes each component of the model stated in
Section~C.5.

\subsubsection*{C.6.1 Objective Function}

The objective~\eqref{eq:pricing_obj} maximizes a single user-selected
KPI evaluated over the PepsiCo product set $\mathcal{P}$. The candidate
KPIs are revenue~\eqref{eq:revenue}, profit~\eqref{eq:profit},
volume~\eqref{eq:total_volume}, margin~\eqref{eq:margin_pricing}, market
share~\eqref{eq:market_share_pricing}, and net revenue per unit
weight~\eqref{eq:nrp}. The nonlinear and non-convex nature of these
functions, arising from the exponential volume model, the logarithmic
elasticity terms, and the PPP discontinuities, motivates the use of a
metaheuristic solver.

\subsubsection*{C.6.2 Decision Variables and Bounds}

The decision variable $x_g$ has a market-dependent interpretation.
In percentage mode (e.g., the U.S. market), $x_g \in \mathbb{R}$ is a
multiplicative factor on the baseline price, as in~\eqref{eq:price_pct}.
In absolute mode (e.g., the Mexico market), $x_g \in \mathbb{Z}$ is the
new shelf price expressed directly in currency units, as
in~\eqref{eq:price_abs}. In both cases, $x_g$ is bounded by user-supplied
limits via~\eqref{eq:var_bounds}.

\subsubsection*{C.6.3 Price Computation}

Given a decision vector $\mathbf{x}$, final shelf prices for all products
are computed in three steps. First, PepsiCo new prices are obtained
from~\eqref{eq:price_pct} or~\eqref{eq:price_abs}, and the corresponding
percentage change is computed via~\eqref{eq:pep_pctchange}. Second,
competitor prices adjust in proportion to their associated PepsiCo
``hero'' product through~\eqref{eq:comp_follow}. Finally,
\eqref{eq:final_price} produces the final shelf prices used in the
volume function and downstream KPIs.

\subsubsection*{C.6.4 Volume Function}

The predicted volume in~\eqref{eq:volume} depends on two components:
a cross-price elasticity term~\eqref{eq:elast_component} and a
psychological pricing term~\eqref{eq:ppp_component}. The exponential
structure ensures that at baseline prices
($\mathbf{p}^{\mathrm{f}} = \mathbf{p}^{0}$), both terms vanish and
$V_i = V_i^{0}$.

In~\eqref{eq:elast_component}, $E_{ij}$ is the $(i,j)$-th entry of the
elasticity matrix $\mathbf{E}$. The diagonal entries $E_{ii}$ capture
own-price elasticity (typically negative), while off-diagonal entries
$E_{ij}$, $i \neq j$, capture cross-price effects (positive for
substitutes, negative for complements).

The PPP function~\eqref{eq:ppp} captures discrete demand shifts that
occur when a product's price crosses consumer-sensitive thresholds
(e.g., \$1.99, 20~Pesos), where $\unrhd$ denotes a market-specific
comparison operator ($\geq$ or $>$), and the sum aggregates the PPP
coefficients $\psi_{i,k}$ for all thresholds met by the current price.
Because the PPP function is a step function of prices, it introduces
discontinuities into the volume response surface.

\subsubsection*{C.6.5 Financial KPIs}

Equations~\eqref{eq:units}--\eqref{eq:nrp} define the financial KPIs
over a product subset $\mathcal{S} \subseteq \mathcal{A}$ (typically
$\mathcal{S} = \mathcal{P}$). Sales units are obtained from predicted
volumes via pack-size adjustment~\eqref{eq:units}. The sell-in price
is updated~\eqref{eq:sip_update}, capturing the manufacturer-to-retailer
price pass-through. Revenue~\eqref{eq:revenue}, profit~\eqref{eq:profit},
margin~\eqref{eq:margin_pricing}, total volume~\eqref{eq:total_volume},
market share~\eqref{eq:market_share_pricing}, and net revenue per unit
weight~\eqref{eq:nrp} follow from these primitives.

\subsubsection*{C.6.6 Financial KPI Bounds}

Constraint~\eqref{eq:kpi_bounds} bounds each financial KPI relative to
its baseline value. Users configure $\underline{\alpha}^{K}$ and
$\overline{\alpha}^{K}$ through the interface for any combination of
revenue~\eqref{eq:revenue}, profit~\eqref{eq:profit},
volume~\eqref{eq:total_volume}, margin~\eqref{eq:margin_pricing},
market share~\eqref{eq:market_share_pricing}, and net revenue per unit
weight~\eqref{eq:nrp}. For instance, setting $\underline{\alpha}^{R} = 0.95$
ensures that revenue does not decrease by more than 5\% from the baseline.

\subsubsection*{C.6.7 Price Ladder Constraint}

The ladder constraint~\eqref{eq:ladder} enforces a
price-per-unit-weight hierarchy within each brand unit: larger pack
sizes must offer a lower (or equal) price per unit weight than smaller
sizes. The factor $\gamma \leq 1$ is a configurable tolerance
(default $\gamma = 1$, i.e., strict ordering).

\subsubsection*{C.6.8 Tier Pricing Constraint}

The tier constraint~\eqref{eq:tier} enforces price hierarchies between
product tiers within the portfolio, ensuring that premium-tier products
are priced at or above lower-tier products.

\subsubsection*{C.6.9 Pricing Line Constraint}

Constraint~\eqref{eq:pricing_line} requires products within the same
pricing line to receive equal absolute price increments, ensuring
uniform price movements across products that share a commercial
pricing group (e.g., the same brand family at different retailers).

\subsubsection*{C.6.10 Average Price Change Constraint}

Constraint~\eqref{eq:avg_price} bounds the average percentage price
change across the PepsiCo portfolio.

\subsubsection*{C.6.11 Volume-Weighted Average Price Change Constraint}

The volume-weighted variant of~\eqref{eq:avg_price}, defined
through~\eqref{eq:wtd_price} and bounded by~\eqref{eq:wtd_price_bound},
uses the \emph{new} predicted volumes $V_i(\mathbf{p}^{\mathrm{f}})$ as
weights. This constraint is endogenous: the weights depend on the
decision variables through the nonlinear volume function~\eqref{eq:volume},
adding further nonlinearity to the feasible region.

%
%



\bibliographystyle{informs2014} 


\begin{thebibliography}{24}
\providecommand{\natexlab}[1]{#1}
\providecommand{\url}[1]{\texttt{#1}}
\providecommand{\urlprefix}{URL }

\bibitem[{Bertsimas \protect\BIBand{} Kallus(2020)}]{bertsimas2020predictive}
Bertsimas D, Kallus N (2020) From predictive to prescriptive analytics. \emph{Management Science} 66(3):1025--1044.

\bibitem[{Bertsimas \protect\BIBand{} Tsitsiklis(1997)}]{bertsimas1997introduction}
Bertsimas D, Tsitsiklis JN (1997) \emph{Introduction to linear optimization}, volume~6 (Athena Scientific Belmont, MA).

\bibitem[{Carpenter et~al.(2017)Carpenter, Gelman, Hoffman, Lee, Goodrich, Betancourt, Brubaker, Guo, Li, \protect\BIBand{} Riddell}]{carpenter2017stan}
Carpenter B, Gelman A, Hoffman MD, Lee D, Goodrich B, Betancourt M, Brubaker M, Guo J, Li P, Riddell A (2017) Stan: A probabilistic programming language. \emph{Journal of Statistical Software} 76(1):1--32.

\bibitem[{Cohen et~al.(2018)Cohen, Lobel, \protect\BIBand{} Perakis}]{cohen2018dynamic}
Cohen MC, Lobel R, Perakis G (2018) Dynamic pricing through data sampling. \emph{Production and Operations Management} 27(6):1074--1088.

\bibitem[{Deng et~al.(2023)Deng, Zhang, Wang, Wang, Zhang, Wang, Zhao, Qi, Yang, \protect\BIBand{} Peng}]{deng2023alibaba}
Deng Y, Zhang X, Wang T, Wang L, Zhang Y, Wang X, Zhao S, Qi Y, Yang G, Peng X (2023) Alibaba realizes millions in cost savings through integrated demand forecasting, inventory management, price optimization, and product recommendations. \emph{INFORMS Journal on Applied Analytics} 53(1):32--46.

\bibitem[{Dimitrov \protect\BIBand{} Durai(2023)}]{digital_cpg}
Dimitrov S, Durai A (2023) Digitizing revenue growth management for scale and impact. \emph{MathCo} \urlprefix\url{https://mathco.com/article/digitizing-rgm-for-scale-and-impact/}.

\bibitem[{Fattahi et~al.(2022)Fattahi, Li, \protect\BIBand{} Sahin}]{fattahi2022customer}
Fattahi A, Li Y, Sahin O (2022) Customer-driven bundle promotion optimization at scale. \emph{Johns Hopkins Carey Business School Research Paper} (22-14).

\bibitem[{Ferreira et~al.(2018)Ferreira, Simchi-Levi, \protect\BIBand{} Wang}]{ferreira2018online}
Ferreira KJ, Simchi-Levi D, Wang H (2018) Online network revenue management using thompson sampling. \emph{Operations Research} 66(6):1586--1602.

\bibitem[{Friedman(2001)}]{friedman2001greedy}
Friedman JH (2001) Greedy function approximation: a gradient boosting machine. \emph{Annals Of Statistics} 1189--1232.

\bibitem[{Gelman et~al.(1995)Gelman, Carlin, Stern, \protect\BIBand{} Rubin}]{gelman1995bayesian}
Gelman A, Carlin JB, Stern HS, Rubin DB (1995) \emph{Bayesian data analysis} (Chapman and Hall/CRC).

\bibitem[{{Gurobi Optimization, LLC}(2024)}]{gurobi}
{Gurobi Optimization, LLC} (2024) {Gurobi Optimizer Reference Manual}. \urlprefix\url{https://www.gurobi.com}.

\bibitem[{Hazan et~al.(2021)Hazan, Br{\'e}g{\'e}, Verwaerde, \protect\BIBand{} Bassoulet}]{bcg2021aipricing}
Hazan J, Br{\'e}g{\'e} C, Verwaerde JS, Bassoulet A (2021) Why ai transformations should start with pricing. \emph{Boston Consulting Group} \urlprefix\url{https://www.bcg.com/publications/2021/ai-pricing-tranformations}.

\bibitem[{Hormby et~al.(2010)Hormby, Morrison, Dave, Meyers, \protect\BIBand{} Tenca}]{hormby2010marriott}
Hormby S, Morrison J, Dave P, Meyers M, Tenca T (2010) Marriott international increases revenue by implementing a group pricing optimizer. \emph{Interfaces} 40(1):47--57.

\bibitem[{Ito \protect\BIBand{} Fujimaki(2017)}]{ito2017optimization}
Ito S, Fujimaki R (2017) Optimization beyond prediction: Prescriptive price optimization. \emph{Proceedings of the 23rd ACM SIGKDD international conference on knowledge discovery and data mining}, 1833--1841.

\bibitem[{Ke et~al.(2017)Ke, Meng, Finley, Wang, Chen, Ma, Ye, \protect\BIBand{} Liu}]{ke2017lightgbm}
Ke G, Meng Q, Finley T, Wang T, Chen W, Ma W, Ye Q, Liu TY (2017) Lightgbm: A highly efficient gradient boosting decision tree. \emph{Advances in Neural Information Processing Systems}, 3146--3154.

\bibitem[{Kucukelbir et~al.(2017)Kucukelbir, Tran, Ranganath, Gelman, \protect\BIBand{} Blei}]{kucukelbir2017advi}
Kucukelbir A, Tran D, Ranganath R, Gelman A, Blei DM (2017) Automatic differentiation variational inference. \emph{Journal of Machine Learning Research} 18(14):1--45.

\bibitem[{Ma et~al.(2021)Ma, Simchi-Levi, \protect\BIBand{} Zhao}]{ma2021dynamic}
Ma W, Simchi-Levi D, Zhao J (2021) Dynamic pricing (and assortment) under a static calendar. \emph{Management Science} 67(4):2292--2313.

\bibitem[{Miao(2020)}]{miao2020data}
Miao S (2020) \emph{Data-Driven Optimization in Revenue Management: Pricing, Assortment Planning, and Demand Learning}. Ph.D. thesis.

\bibitem[{PepsiCo(2025)}]{PepsiCo}
PepsiCo (2025) \urlprefix\url{https://www.pepsico.com/who-we-are/about-pepsico}.

\bibitem[{Storn \protect\BIBand{} Price(1997)}]{storn1997differential}
Storn R, Price K (1997) Differential evolution--a simple and efficient heuristic for global optimization over continuous spaces. \emph{Journal of Global Optimization} 11(4):341--359.

\bibitem[{Talluri \protect\BIBand{} Van~Ryzin(2006)}]{talluri2006theory}
Talluri KT, Van~Ryzin GJ (2006) \emph{The theory and practice of revenue management}, volume~68 (Springer Science \& Business Media).

\bibitem[{Virtanen et~al.(2020)Virtanen, Gommers, Oliphant, Haberland, Reddy, Cournapeau, Burovski, Peterson, Weckesser, Bright, {van der Walt}, Brett, Wilson, Millman, Mayorov, Nelson, Jones, Kern, Larson, Carey, Polat, Feng, Moore, {VanderPlas}, Laxalde, Perktold, Cimrman, Henriksen, Quintero, Harris, Archibald, Ribeiro, Pedregosa, {van Mulbregt}, \protect\BIBand{} {SciPy 1.0 Contributors}}]{2020SciPy-NMeth}
Virtanen P, Gommers R, Oliphant TE, Haberland M, Reddy T, Cournapeau D, Burovski E, Peterson P, Weckesser W, Bright J, {van der Walt} SJ, Brett M, Wilson J, Millman KJ, Mayorov N, Nelson ARJ, Jones E, Kern R, Larson E, Carey CJ, Polat {\.I}, Feng Y, Moore EW, {VanderPlas} J, Laxalde D, Perktold J, Cimrman R, Henriksen I, Quintero EA, Harris CR, Archibald AM, Ribeiro AH, Pedregosa F, {van Mulbregt} P, {SciPy 10 Contributors} (2020) {{SciPy} 1.0: Fundamental Algorithms for Scientific Computing in Python}. \emph{Nature Methods} 17:261--272, \urlprefix\url{http://dx.doi.org/10.1038/s41592-019-0686-2}.

\bibitem[{Wicaksono \protect\BIBand{} Karimi(2008)}]{wicaksono2008piecewise}
Wicaksono DS, Karimi IA (2008) Piecewise MILP under- and overestimators for global optimization of bilinear programs. \emph{AIChE Journal} 54(4):991--1008.

\bibitem[{Wolsey \protect\BIBand{} Nemhauser(1999)}]{wolsey1999integer}
Wolsey LA, Nemhauser GL (1999) \emph{Integer and combinatorial optimization} (John Wiley \& Sons).

\end{thebibliography}





\end{document}